\documentclass [reqno]{amsart}
\usepackage {mathptmx}
\usepackage {amsmath,amssymb}
\usepackage {color,graphicx}
\usepackage [all]{xy}

\newcommand {\CC}{{\mathbb C}}

\newcommand {\PP}{{\mathbb P}}
\newcommand {\QQ}{{\mathbb Q}}
\newcommand {\RR}{{\mathbb R}}
\newcommand {\SM}{{\mathbb S}}
\newcommand {\ZZ}{{\mathbb Z}}

\DeclareMathOperator {\ft}{ft}
\DeclareMathOperator {\ev}{ev}
\DeclareMathOperator {\id}{id}

\DeclareMathOperator {\val}{val}
\DeclareMathOperator {\dist}{dist}

\DeclareMathOperator {\im}{Im}
\DeclareMathOperator {\Trop}{Trop}
\DeclareMathOperator {\ini}{in}
\DeclareMathOperator {\mult}{mult}
\DeclareMathOperator {\trop}{trop}
\DeclareMathOperator {\lab}{lab}
\DeclareMathOperator {\pr}{pr}

\DeclareSymbolFont {mysymbols}{OMS}{cmsy}{m}{n}
\DeclareMathSymbol {\calI}{\mathalpha}{mysymbols}{`I}

\DeclareSymbolFont {mylargesymbols}{OMX}{cmex}{m}{n}
\DeclareMathSymbol {\dunion}{\mathop}{mylargesymbols}{"60}

\newcommand {\calM}{{\mathcal M}}
\newcommand {\df}[1]{\textsl {#1}}

\newcommand {\preprint}[2]{preprint \discretionary {#1/}{#2}{#1/#2}}

\renewenvironment {enumerate}%
  {\rule{1mm}{0mm}\begin {oldenumerate}%
    \parskip1ex plus0.5ex \itemsep 0mm \parindent 0mm}%
  {\end {oldenumerate}}

\renewenvironment {itemize}%
  {\rule{1mm}{0mm}\begin {olditemize}%
    \parskip1ex plus0.5ex \itemsep 0mm \parindent 0mm}%
  {\end {olditemize}}

\theoremstyle {plain}
\newtheorem {theorem}{Theorem}[section]
\newtheorem {proposition}[theorem]{Proposition}
\newtheorem {lemma}[theorem]{Lemma}
\newtheorem {corollary}[theorem]{Corollary}
\theoremstyle {definition}
\newtheorem {definition}[theorem]{Definition}
\theoremstyle {remark}
\newtheorem {remark}[theorem]{Remark}
\newtheorem {remdef}[theorem]{Remark and Definition}
\newtheorem {example}[theorem]{Example}
\newtheorem {construction}[theorem]{Construction}

\begin {document}

  \title [Tropical fans and the moduli spaces of tropical curves]{Tropical
  fans and the moduli spaces of tropical curves}
\author {Andreas Gathmann, Michael Kerber, and Hannah Markwig}

\address {Andreas Gathmann, Fachbereich Mathematik, TU Kaiserslautern, Postfach
  3049, 67653 Kaiserslautern, Germany}
\email {andreas@mathematik.uni-kl.de}
\address {Michael Kerber, Fachbereich Mathematik, TU Kaiserslautern, Postfach
  3049, 67653 Kaiserslautern, Germany}
\email {mkerber@mathematik.uni-kl.de}
\address {Hannah Markwig, Institute for Mathematics and its Applications (IMA),
  University of Minne\-sota, Lind Hall 400, 207 Church Street SE, Minneapolis,
  MN 55455, USA}
\email {markwig@ima.umn.edu}

\begin {abstract}
  We give a rigorous definition of tropical fans (the ``local building blocks
  for tropical varieties'') and their morphisms. For a morphism of
  tropical fans of the same dimension we show that the number of inverse images
  (counted with suitable tropical multiplicities) of a point in the target does
  not depend on the chosen point --- a statement that can be viewed as one of
  the important first steps of tropical intersection theory. As an application
  we consider the moduli spaces of rational tropical curves (both abstract and
  in some $ \RR^r $) together with the evaluation and forgetful morphisms.
  Using our results this gives new, easy, and unified proofs of various
  tropical independence statements, e.g.\ of the fact that the numbers of
  rational tropical curves (in any $ \RR^r $) through given points are
  independent of the points.
\end {abstract}

\maketitle

  \section {Introduction}

Tropical geometry has been used recently for many applications in enumerative
geometry, both to reprove previously known results (e.g.\ the Caporaso-Harris
formula \cite {GM2} and Kontsevich's formula \cite {GM} for plane curves) and
to find new methods for counting curves (e.g.\ Mikhalkin's lattice path
algorithm \cite {M0}). All these applications involve counting curves with
certain conditions, and at some point require an argument why the resulting
numbers remain invariant under deformations of the conditions --- e.g.\ why the
numbers of tropical plane rational curves of degree $d$ through $ 3d-1 $ points
do not depend on the chosen points (when counted with the appropriate tropical
multiplicities).

In classical algebraic geometry, the corresponding independence statements
simply follow from the fact that the enumerative numbers can be interpreted as
intersection numbers of cycles on suitable moduli spaces, and that intersection
numbers are always invariant under deformations. Tropically however there is no
well-established intersection theory yet, and thus so far it was necessary to
find alternative proofs of the independence statements in each case --- either
by tedious ad-hoc computations (as e.g.\ in \cite {GM}) or by relating the
tropical numbers to classical ones for which the invariance is known (as e.g.\
in \cite {M0}).

The goal of our paper is to fill this gap by establishing the very basics of a
``tropical intersection theory'' up to the point needed so that at least some
of the above independence statements follow from general principles of this
theory, just as in classical algebraic geometry.

To do so we start in chapter \ref {sec-tropfan} by giving a rigorous definition
of the notion of tropical fans (roughly speaking tropical varieties in some
vector space all of whose cells are cones with apex at the origin, i.e.\ a
``local picture'' of a general tropical variety), morphisms between them, and
the image fan of a morphism. Our most important result in this chapter is
corollary \ref {cor-image} which states that for a morphism of tropical fans of
the same dimension, with the target being irreducible, the sum of the
multiplicities of the inverse images of a general point $P$ in the target is
independent of the choice of $P$. In chapter \ref {sec-m0n} we establish the
moduli spaces of abstract $n$-marked rational tropical curves as tropical fans,
and show that the forgetful maps between them are morphisms of fans. We then
use these results in chapter \ref {sec-smap} to construct tropical analogues of
the spaces of rational stable maps to a toric variety. Again, these spaces will
be tropical fans, with the evaluation and forgetful maps being morphisms of
fans. In chapter \ref {sec-cor} finally, we present two examples of our
results by giving new, easy, and generalized proofs of two statements in
tropical enumerative geometry that have already occurred earlier in the
literature: the statement that the numbers of rational tropical curves in some
$ \RR^r $ of given degree and through fixed affine linear subspaces in general
position do not depend on the position of the subspaces, and the statement that
the degree of the combined evaluation and forgetful map occurring in the proof
of the tropical Kontsevich formula is independent of the choice of the points
(see \cite {GM} proposition 4.4). Whereas earlier tropical proofs of these
statements have been very complicated, our new proofs are now just easy
applications of our general statement in corollary \ref {cor-image}.

We believe that our work results not only in a better understanding of the
moduli spaces of tropical curves, but also of tropical varieties and their
intersection-theoretical properties in general (see \cite{M}). Subsequent work on a tropical
intersection theory building up on the principles of this paper is in progress.

We should mention that a theorem very similar to our corollary \ref {cor-image}
has recently been proven independently by Bernd Sturmfels and Jenia Tevelev
(see \cite {ST} theorem 1.1). In their paper the authors restrict to tropical
fans that are tropicalizations of algebraic varieties (note that these objects
are just called ``tropicalizations'' in their work, whereas their term
``tropical fan'' has a meaning different from ours). As a result, they are able
not only to prove the independence statement of our corollary \ref {cor-image}
but also to show that this invariant is equal to the degree of the
corresponding morphism of algebraic varieties. On the other hand, our result is
of course applicable in more generality as it does not need the (rather strong)
requirement on the tropical fans to be tropicalizations of algebraic varieties.
Another consequence of the different point of view in \cite {ST} compared to
our work is that the proofs are entirely different: whereas Sturmfels and
Tevelev prove their theorem 1.1 by relating the tropical set-up to the
algebraic one our proof of corollary \ref {cor-image} is entirely combinatorial
and does not use any algebraic geometry.

The second author would like to thank the Institute for Mathematics and its
Applications (IMA) in Minneapolis for hospitality. The third author would like
to thank Ionut Ciocan-Fontanine, Diane Maclagan, and Grisha Mikhalkin for many
helpful discussions.

  \section {Tropical fans} \label {sec-tropfan}

Throughout this section $ \Lambda $ will denote a finitely generated free
abelian group (i.e.\ a group isomorphic to $ \ZZ^N $ for some $ N \ge 0 $) and
$ V := \Lambda \otimes_\ZZ \RR $ the corresponding real vector space, so that $
\Lambda $ can be considered as a lattice in $V$. The dual lattice in the dual
vector space  $V^\vee $ will be denoted $ \Lambda^\vee $.

\begin {definition}[Cones] \label {def-cone}
  A \df {cone} is a subset $ \sigma \subset V $ given by finitely many linear
  integral equalities and (non-strict) inequalities, i.e.\ a set of the form
    \[ \sigma = \{ x \in V ;\;
         \mbox {$ f_i (x) = 0 $ for all $ i=1,\dots,n $ and
                $ g_j (x) \ge 0 $ for all $ j=1,\dots,m $} \} \tag {$*$} \]
  for some $ f_1,\dots,f_n,g_1,\dots,g_m \in \Lambda^\vee \subset V^\vee $.
  We denote by $ V_\sigma \subset V $ the smallest vector subspace of $V$ that
  contains $ \sigma $, and by $ \Lambda_\sigma := V_\sigma \cap \Lambda \subset
  \Lambda $ the smallest sublattice of $ \Lambda $ that contains $ \sigma \cap
  \Lambda $. The \df {dimension} of $ \sigma $ is defined to be $ \dim \sigma
  := \dim V_\sigma $.

  Let $ \sigma $ be a cone in $V$. A cone $ \tau \subset \sigma $ that
  can be
  obtained from $ \sigma $ by changing some (maybe none) of the inequalities
  $ g_j (x) \ge 0 $ in $ (*) $ to equalities $ g_j(x) = 0 $ is called a \df
  {face} of $ \sigma $. We write this as $ \tau \le \sigma $ (or as $ \tau <
  \sigma $ if in addition $ \tau \subsetneq \sigma $). We obviously have $
  V_\tau \subset V_\sigma $ and $ \Lambda_\tau \subset \Lambda_\sigma $ in this
  case.
\end {definition}

\begin {remdef} \label {remdef-cone}
  It is well-known that a cone can be described equivalently as a subset of $V$
  that can be written as
    \[ \sigma = \{ \lambda_1 u_1 + \cdots + \lambda_n u_n ;\;
         \lambda_1,\dots,\lambda_n \in \RR_{\ge 0} \} \]
  for some integral vectors $ u_1,\dots,u_n \in \Lambda $. In this case we say
  that the cone $ \sigma $ is \df {generated} by $ u_1,\dots,u_n $. A cone is
  called an \df {edge} if it can be generated by one vector, and \df
  {simplicial} if it can be generated by $ \dim \sigma $ vectors (that are then
  necessarily linearly independent).

  If $ \sigma $ is generated by $ u_1,\dots,u_n \in \Lambda $ then each face of
  $ \sigma $ can be generated by a suitable subset of $ \{ u_1,\dots,u_n \} $.
  In fact, if $ \sigma $ is simplicial then there is a 1:1 correspondence
  between faces of $ \sigma $ and subsets of $ \{ u_1,\dots,u_n \} $. In
  particular, such a simplicial cone has among its faces exactly $n$ edges,
  namely $ \RR_{\ge 0} \cdot u_1,\dots,\RR_{\ge 0} \cdot u_n $.
\end {remdef}

\begin {construction}[Normal vectors] \label {constr-cone}
  Let $ \tau < \sigma $ be cones in $V$ with $ \dim \tau = \dim \sigma - 1 $.
  By definition there is a linear form $ g \in \Lambda^\vee $ that is zero on
  $ \tau $, non-negative on $ \sigma $, and not identically zero on $ \sigma $.
  Then $g$ induces an isomorphism $ V_\sigma / V_\tau \to \RR $ that is
  non-negative and not identically zero on $ \sigma/ V_\tau  $, showing that the cone $
  \sigma / V_\tau $ lies in a unique half-space of $ V_\sigma / V_\tau \cong \RR $. As
  moreover $ \Lambda_\sigma / \Lambda_\tau \subset V_\sigma / V_\tau $ is
  isomorphic to $ \ZZ $ there is a unique generator of $ \Lambda_\sigma /
  \Lambda_\tau $ lying in this half-space. We denote it by $ u_{\sigma/\tau}
  \in \Lambda_\sigma / \Lambda_\tau $ and call it the \df {(primitive) normal
  vector} of $ \sigma $ relative to $ \tau $.
\end {construction}

\begin {definition}[Fans] \label {def-fan}
  A \df {fan} $X$ in $V$ is a finite set of cones in $V$ such that
  \begin {enumerate}
  \item \label {def-fan-a}
    all faces of the cones in $X$ are also in $X$; and
  \item \label {def-fan-b}
    the intersection of any two cones in $X$ is a face of each (and hence by
    \ref {def-fan-a} also in $X$).
  \end {enumerate}
  The set of all $k$-dimensional cones of $X$ will be denoted $ X^{(k)} $. The biggest
  dimension of a cone in $X$ is called the \df {dimension} $ \dim X $ of $X$;
  we say that $X$ is \df {pure-dimensional} if each inclusion-maximal cone in
  $X$ has this dimension. A fan is called \df {simplicial} if all of its cones
  are simplicial. The union of all cones in $X$ will be denoted $ |X| \subset
  V $.
\end {definition}

\begin {example} \label {ex-fan}
  \begin {enumerate}
  \item \label {ex-fan-a}
    In contrast to the fans considered in toric geometry we do not require our
    cones to be \textsl {strictly} convex, i.e.\ a cone might contain a
    straight line through the origin. So e.g.\ $ \{V\} $ is a fan in our
    sense with only one cone (that has no faces except itself). By abuse of
    notation we will denote it simply by $V$.
  \item \label {ex-fan-b}
    For $ f \in \Lambda^\vee \backslash \{0\} $ the three cones
      \[ \{ x \in V;\; f(x) = 0 \}, \quad
         \{ x \in V;\; f(x) \ge 0 \}, \quad
         \{ x \in V;\; f(x) \le 0 \} \]
    form a fan. We call it the \df {half-space fan} $ H_f $.
  \item \label {ex-fan-c}
    Let $ \Lambda = \ZZ^n $ (and thus $ V = \RR^n $), let $ u_1,\dots,u_n $ be
    a basis of $ \Lambda $, and set $ u_0 := -u_1-\cdots-u_n $. For each
    subset $ I \subsetneq \{ 0,\dots,n \} $ we denote by $ \sigma_I $ the
    simplicial cone of dimension $ |I| $ in $V$ spanned by the vectors $ u_i $
    for $ i \in I $. Now fix an integer $ k \in \{0,\dots,n\} $ and let $ L_k^n
    $ be the set of all cones $ \sigma_I $ for $ |I| \le k $. Then the
    conditions of definition \ref {def-fan} are satisfied since the faces of $
    \sigma_I $ are precisely the $ \sigma_J $ with $ J \subset I $ by remark
    \ref {remdef-cone}, and moreover we have $ \sigma_{I_1} \cap \sigma_{I_2} =
    \sigma_{I_1 \cap I_2} $ for all $ I_1,I_2 \subsetneq \{ 0,\dots,n \} $.
    Hence $ L_k^n $ is a fan. It is clearly of pure dimension $k$. For $ n=2 $
    we obtain the following picture:
    \begin {center} \begin{picture}(0,0)%
\includegraphics{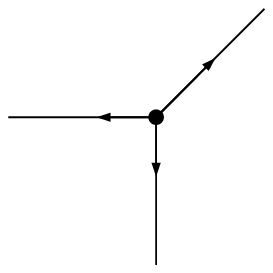}%
\end{picture}%
\setlength{\unitlength}{4144sp}%
\begingroup\makeatletter\ifx\SetFigFontNFSS\undefined%
\gdef\SetFigFontNFSS#1#2#3#4#5{%
  \reset@font\fontsize{#1}{#2pt}%
  \fontfamily{#3}\fontseries{#4}\fontshape{#5}%
  \selectfont}%
\fi\endgroup%
\begin{picture}(1254,1369)(412,-833)
\put(427,-124){\makebox(0,0)[rb]{\smash{{\SetFigFontNFSS{10}{12.0}{\familydefault}{\mddefault}{\updefault}{\color[rgb]{0,0,0}$\sigma_{\{1\}}$}%
}}}}
\put(1157,-769){\makebox(0,0)[lb]{\smash{{\SetFigFontNFSS{10}{12.0}{\familydefault}{\mddefault}{\updefault}{\color[rgb]{0,0,0}$\sigma_{\{2\}}$}%
}}}}
\put(892,  3){\makebox(0,0)[b]{\smash{{\SetFigFontNFSS{10}{12.0}{\familydefault}{\mddefault}{\updefault}{\color[rgb]{0,0,0}$u_1$}%
}}}}
\put(1098,-353){\makebox(0,0)[rb]{\smash{{\SetFigFontNFSS{10}{12.0}{\familydefault}{\mddefault}{\updefault}{\color[rgb]{0,0,0}$u_2$}%
}}}}
\put(1402, 78){\makebox(0,0)[lb]{\smash{{\SetFigFontNFSS{10}{12.0}{\familydefault}{\mddefault}{\updefault}{\color[rgb]{0,0,0}$u_0$}%
}}}}
\put(1172,-161){\makebox(0,0)[lb]{\smash{{\SetFigFontNFSS{10}{12.0}{\familydefault}{\mddefault}{\updefault}{\color[rgb]{0,0,0}$\sigma_\emptyset$}%
}}}}
\put(1651,374){\makebox(0,0)[lb]{\smash{{\SetFigFontNFSS{10}{12.0}{\familydefault}{\mddefault}{\updefault}{\color[rgb]{0,0,0}$\sigma_{\{0\}}$}%
}}}}
\put(1621,-376){\makebox(0,0)[b]{\smash{{\SetFigFontNFSS{10}{12.0}{\familydefault}{\mddefault}{\updefault}{\color[rgb]{0,0,0}$\sigma_{\{0,2\}}$}%
}}}}
\put(856,389){\makebox(0,0)[b]{\smash{{\SetFigFontNFSS{10}{12.0}{\familydefault}{\mddefault}{\updefault}{\color[rgb]{0,0,0}$\sigma_{\{0,1\}}$}%
}}}}
\put(676,-556){\makebox(0,0)[b]{\smash{{\SetFigFontNFSS{10}{12.0}{\familydefault}{\mddefault}{\updefault}{\color[rgb]{0,0,0}$\sigma_{\{1,2\}}$}%
}}}}
\end{picture}%
 \end {center}
    where $ L_0^2 $ is just the origin, $ L_1^2 $ is the ``tropical line''
    given by the origin together with the three edges spanned by $ u_0 $, $ u_1
    $, $ u_2 $, and $ |L_2^2|=\RR^2 $.
  \item \label {ex-fan-d}
    If $X$ and $X'$ are two fans in $ V = \Lambda \otimes \RR $ and $ V' =
    \Lambda' \otimes \RR $, respectively, it is checked immediately that their
    product $ \{ \sigma \times \sigma' ;\; \sigma \in X ,\; \sigma' \in X' \} $
    is a fan in $ V \times V' $. We call it the \df {product} of the two fans
    and denote it by $ X \times X' $. We obviously have $ |X \times X'| = |X|
    \times |X'| $.
  \item \label {ex-fan-e}
    In the same way, if $X$ and $X'$ are two fans in $V$ it is checked
    immediately that their intersection $ \{ \sigma \cap \sigma' ;\; \sigma \in
    X ,\; \sigma' \in X' \} $ is also a fan in $V$. We call it the \df
    {intersection} of the two fans and denote it by $ X \cap X' $. It is clear
    that $ |X \cap X'| = |X| \cap |X'| $.
  \end {enumerate}
\end {example}

\begin {definition}[Subfans] \label {def-subfan}
  Let $X$ be a fan in $V$. A fan $Y$ in $V$ is called a \df {subfan} of $X$
  (denoted $ Y \subset X $) if each cone of $Y$ is contained in a cone of $X$.
  In this case we denote by $ C_{Y,X}: Y \to X $ the map that sends a cone $
  \sigma \in Y $ to the (unique) inclusion-minimal cone of $X$ that contains $
  \sigma $. Note that for a subfan $ Y \subset X $ we obviously have $ |Y|
  \subset |X| $ and $ \dim C_{Y,X} (\sigma) \ge \dim \sigma $ for all $ \sigma
  \in Y $.
\end {definition}

\begin {example} \label {ex-subfan}
  The intersection of two fans as in example \ref {ex-fan} \ref {ex-fan-e} is
  clearly a subfan of both.
\end {example}

\begin {definition}[Weighted and tropical fans] \label {def-tropfan}
  A \df {weighted fan} in $V$ is a pair $ (X,\omega_X) $ where $X$ is a fan of
  some pure dimension $n$ in $V$, and $ \omega_X: X^{(n)} \to \ZZ_{>0} $ is a
  map. We call $ \omega_X(\sigma) $ the \df {weight} of the cone $ \sigma \in
  X^{(n)} $ and write it simply as $ \omega(\sigma) $ if no confusion can
  result. Also, by abuse of notation we will write a weighted fan $
  (X,\omega_X) $ simply as $X$ if the weight function $ \omega_X $ is clear
  from the context.

  A \df {tropical fan} in $V$ is a weighted fan $ (X,\omega_X) $ in $V$ such
  that for all $ \tau \in X^{(\dim X-1)} $ the \df {balancing condition}
    \[ \sum_{\sigma > \tau} \omega_X(\sigma) \cdot u_{\sigma/\tau} = 0
         \quad \in V / V_\tau \]
  holds (where $ u_{\sigma/\tau} $ denotes the primitive normal vector of
  construction \ref {constr-cone}).
\end {definition}

\begin {example} \label {ex-tropfan}
  \begin {enumerate}
  \item \label {ex-tropfan-a}
    Of course, the fan $V$ of example \ref {ex-fan} \ref {ex-fan-a} with the
    weight function $ \omega(V) := 1 $ is a tropical fan (the balancing
    condition being trivial in this case). The half-space fans $ H_f $ of
    example \ref {ex-fan} \ref {ex-fan-b} are tropical fans as well if we
    define the weight of both maximal cones to be 1: the balancing condition
    around the cone $ \{ x \in V ;\; f(x) = 0 \} $ holds because the primitive
    normal vectors of the two adjacent maximal cones are negatives of each
    other. We will see in example \ref {ex-marked} that the fans $ L_k^n $ of
    example \ref {ex-fan} \ref {ex-fan-c} are also tropical fans if the weights
    of all maximal cones are defined to be 1.
  \item \label {ex-tropfan-b}
    Let $ (X,\omega_X) $ be a weighted fan of dimension $n$, and let $ \lambda
    \in \QQ_{>0} $ such that $ \lambda \omega_X(\sigma) \in \ZZ_{>0} $ for all
    $ \sigma \in X^{(n)} $. Then $ (X,\lambda \cdot \omega_X) $ is a weighted
    fan as well; we will denote it by $ \lambda \cdot X $. It is obvious that $
    \lambda \cdot X $ is tropical if and only if $X$ is.
  \item \label {ex-tropfan-c}
    Let $ (X,\omega_X) $ and $ (X',\omega_{X'}) $ be two weighted fans of
    dimensions $n$ and $n'$ in $ V = \Lambda \otimes \RR $ and $ V' = \Lambda'
    \otimes \RR $, respectively. We consider their product $ X \times X' $ of
    example \ref {ex-fan} \ref {ex-fan-d} to be a weighted fan by setting
    $ \omega_{X \times X'}(\sigma \times \sigma') := \omega_X(\sigma) \cdot
    \omega_{X'} (\sigma') $ for $ \sigma \in X^{(n)} $ and $ \sigma' \in
    X'{}^{(n')} $. If moreover $ (X,\omega_X) $ and $ (X',\omega_{X'}) $ are
    tropical then so is $ X \times X' $: the cones in $ (X \times
    X')^{(n+n'-1)} $ are of the form $ \sigma \times \sigma' $ for $ \sigma \in
    X^{(n-1)},\, \sigma' \in X'{}^{(n')} $ or $ \sigma \in X^{(n)},\, \sigma'
    \in X'{}^{(n'-1)} $, and the balancing conditions in these cases are easily
    seen to be induced from the ones of $X$ and $X'$, respectively.
  \end {enumerate}
\end {example}

\begin {definition}[Refinements] \label {def-refine}
  Let $ (X,\omega_X) $ and $ (Y,\omega_Y) $ be weighted fans in $V$. We say
  that $ (Y,\omega_Y) $ is a \df {refinement} of $ (X,\omega_X) $ if
  \begin {itemize}
  \item $ Y \subset X $;
  \item $ |Y|=|X| $ (so in particular $ \dim Y = \dim X $); and
  \item $ \omega_Y(\sigma) = \omega_X (C_{Y,X}(\sigma)) $ for all maximal cones
    $ \sigma \in Y $ (where $ C_{Y,X} $ is as in definition \ref {def-subfan}).
  \end {itemize}
  We say that the two weighted fans $ (X,\omega_X) $ and $ (Y,\omega_Y) $ are
  \df {equivalent} (written $ (X,\omega_X) \cong (Y,\omega_Y) $) if they have a
  common refinement.
\end {definition}

\begin {example} \label {ex-refine}
  \begin {enumerate}
  \item \label {ex-refine-a}
    The half-space fans $ H_f $ of example \ref {ex-fan} \ref {ex-fan-b} are
    all refinements of the trivial fan $V$ of example \ref {ex-fan} \ref
    {ex-fan-a} (with the weights all 1 as in example \ref {ex-tropfan} \ref
    {ex-tropfan-a}).
  \item \label {ex-refine-b}
    Let $ (X,\omega_X) $ be a weighted fan in $V$, and let $Y$ be any fan in
    $V$ with $ |Y| \supset |X| $. Then the intersection $ X \cap Y $ (see
    examples \ref {ex-fan} \ref {ex-fan-e} and \ref {ex-subfan}) is a
    refinement of $ (X,\omega_X) $ by setting $ \omega_{X \cap Y}(\sigma) :=
    \omega_X (C_{X \cap Y,X} (\sigma)) $ for all maximal cones $ \sigma \in X
    \cap Y $. In fact, this construction will be our main source for
    refinements in this paper.

    Note that in the special case when $ (X,\omega_X) $ and $ (Y,\omega_Y) $
    are both weighted fans and $ |Y|=|X| $ we can form both intersections $ X
    \cap Y $ and $ Y \cap X $. The underlying non-weighted fans of these two
    intersections are of course always the same, but the weights may differ as
    they are by construction always induced by the first fan.
  \item \label {ex-refine-c}
    The equivalence of weighted fans is in fact an equivalence relation: if $
    X_1 \cong X_2 $ with common refinement $ Y_1 $, and $ X_2 \cong X_3 $ with
    common refinement $ Y_2 $, then $ Y_1 \cap Y_2 $ as in \ref {ex-refine-b}
    is a common refinement of $ X_1 $ and $ X_3 $ (note that $ Y_1 \cap Y_2 =
    Y_2 \cap Y_1 $ in this case as both $ Y_1 $ and $ Y_2 $ are refinements of
    $ X_2 $).
  \item \label {ex-refine-d}
    Let $ (Y,\omega_Y) $ be a refinement of a weighted fan $ (X,\omega_X) $ of
    dimension $n$. We claim that $ (Y,\omega_Y) $ is tropical if and only if $
    (X,\omega_X) $ is. In fact, for $ \tau \in Y^{(n-1)} $ we have the
    following two cases:
    \begin {itemize}
    \item $ \dim C_{Y,X}(\tau)=n $: Then the $n$-dimensional cones in $Y$ that
      have $ \tau $ as a face must be contained in $ C_{Y,X} (\tau) $. It
      follows that there are exactly two of them, with opposite normal vectors
      and the same weight (equal to $ \omega_X (C_{Y,X}(\tau)) $). In
      particular, the balancing condition for $Y$ always holds at such cones.
    \item $ \dim C_{Y,X}(\tau)=n-1 $: Then the map $ \sigma \mapsto C_{Y,X}
      (\sigma) $ gives a 1:1 correspondence between cones $ \sigma \in
      Y^{(n)} $ with $ \sigma > \tau $ and cones $ \sigma' \in X^{(n)} $ with $
      \sigma' > C_{Y,X}(\tau) $. As this correspondence preserves weights and
      normal vectors the balancing condition for $Y$ at $ \tau $ is equivalent
      to that for $X$ at $ C_{Y,X} (\tau) $. As each $ (n-1) $-dimensional cone
      of $X$ occurs as $ C_{Y,X}(\tau) $ for some $ \tau \in Y^{(n-1)} $ this
      means that the balancing conditions for $X$ and $Y$ are equivalent.
    \end {itemize}
    Consequently, the notion of weighted fans being tropical is well-defined on
    equi\-valence classes.
  \end {enumerate}
\end {example}

\begin {definition}[Marked fans] \label {def-marked}
  A \df {marked fan} in $V$ is a pure-dimensional simplicial fan $X$ in $V$
  together with the data of vectors $ v_\sigma \in (\sigma \backslash \{0\})
  \cap \Lambda $ for all $ \sigma \in X^{(1)} $ (i.e.\ $ v_\sigma $ is an
  integral vector generating the edge $ \sigma $).
\end {definition}

\begin {construction} \label {constr-marked}
  Let $X$ be a marked fan of dimension $n$.
  \begin {enumerate}
  \item \label {constr-marked-a}
    Let $ \sigma \in X^{(k)} $ be a $k$-dimensional cone in $X$. As $ \sigma $
    is simplicial by assumption we know by remark \ref {remdef-cone} that there
    are exactly $k$ edges $ \sigma_1,\dots,\sigma_k \in X^{(1)} $ that are
    faces of $ \sigma $, and that the vectors $ v_{\sigma_1},\dots,v_{\sigma_k}
    $ generating these edges are linearly independent. Hence
      \[ \tilde \Lambda_\sigma
           := \ZZ \, v_{\sigma_1} + \cdots + \ZZ \, v_{\sigma_k} \]
    is a sublattice of $ \Lambda_\sigma $ of full rank, and consequently $
    \Lambda_\sigma / \tilde \Lambda_\sigma $ is a finite abelian group.
    We set $ \omega (\sigma) := | \Lambda_\sigma / \tilde \Lambda_\sigma | \in
    \ZZ_{>0} $. In particular, this makes the marked fan $X$ into a weighted
    fan. In this paper marked fans will always be considered to be weighted
    fans in this way.
  \item \label {constr-marked-b}
    Let $ \sigma \in X^{(k)} $ and $ \tau \in X^{(k-1)} $ with $ \sigma>\tau $.
    As in \ref {constr-marked-a} there are then exactly $ k-1 $ edges in $
    X^{(1)} $ that are faces of $ \tau $, and $k$ edges that are faces of $
    \sigma $. There is therefore exactly one edge $ \sigma' \in X^{(1)} $ that
    is a face of $ \sigma $ but not of $ \tau $. The corresponding vector $
    v_{\sigma'} $ will be denoted $ v_{\sigma/\tau} \in \Lambda $; it can
    obviously also be thought of as a ``normal vector'' of $ \sigma $ relative
    to $ \tau $. Note however that in contrast to the normal vector $
    u_{\sigma/\tau} $ of construction \ref {constr-cone} it is defined in $
    \Lambda $ and not just in $ \Lambda/\Lambda_\tau $, and that it need not be
    primitive.
  \end {enumerate}
\end {construction}

\begin {lemma} \label {lem-marked}
  Let $X$ be a marked fan of dimension $n$ (and hence a weighted fan by
  construction \ref {constr-marked} \ref {constr-marked-a}). Then $X$ is a
  tropical fan if and only if for all $ \tau \in X^{(n-1)} $ the balancing
  condition
    \[ \sum_{\sigma>\tau} v_{\sigma/\tau} = 0 \quad \in V/V_\tau \]
  holds (where the vectors $ v_{\sigma/\tau} $ are as in construction \ref
  {constr-marked} \ref {constr-marked-b}).
\end {lemma}

\begin {proof}
  We have to show that the given equations coincide with the balancing
  condition of definition \ref {def-tropfan}. For this it obviously suffices to
  prove that
    \[ \omega (\sigma) \cdot u_{\sigma/\tau}
         = \omega(\tau) \cdot v_{\sigma/\tau}
           \quad \in \Lambda_\sigma/\Lambda_\tau
           \qquad \mbox {(and hence in $ V_\sigma/V_\tau $)} \]
  for all $ \sigma \in X^{(n)} $ and $ \tau \in X^{(n-1)} $ with $ \sigma>\tau
  $. Using the isomorphism $ \Lambda_\sigma / \Lambda_\tau \cong \ZZ $ of
  construction \ref {constr-cone} this is equivalent to the equation
    \[ \omega (\sigma) = \omega (\tau) \cdot
         | \Lambda_\sigma / (\Lambda_\tau + \ZZ \, v_{\sigma/\tau}) | \]
  in $ \ZZ $, i.e.\ using construction \ref {constr-marked} \ref
  {constr-marked-a} and the relation $ \tilde \Lambda_\sigma = \tilde
  \Lambda_\tau + \ZZ \, v_{\sigma/\tau} $ to the equation
    \[ | \Lambda_\sigma / (\tilde \Lambda_\tau + \ZZ \, v_{\sigma/\tau}) |
         = | \Lambda_\tau / \tilde \Lambda_\tau |
         \cdot | \Lambda_\sigma / (\Lambda_\tau + \ZZ \, v_{\sigma/\tau}) |. \]
  But this follows immediately from the exact sequence
    \[ 0 \;\longrightarrow\; \Lambda_\tau / \tilde \Lambda_\tau
         \;\longrightarrow\; \Lambda_\sigma /
             (\tilde \Lambda_\tau + \ZZ \, v_{\sigma/\tau})
         \;\longrightarrow\; \Lambda_\sigma /
             (\Lambda_\tau + \ZZ \, v_{\sigma/\tau})
         \;\longrightarrow\; 0. \]
\end {proof}

\begin {example} \label {ex-marked}
  We consider again the fans $ L_k^n $ of example \ref {ex-fan} \ref {ex-fan-c}
  and make them into marked fans by setting $ v_{\sigma_{\{i\}}} := u_i $ for $
  i=0,\dots,n $. We claim that $ L_k^n $ becomes a tropical fan in this way. In
  fact, to prove this we just have to check the balancing condition of lemma
  \ref {lem-marked}. By symmetry it suffices to do this at the $ (k-1)
  $-dimensional cone $ \sigma_{\{0,\dots,k-2\}} $. The $k$-dimensional cones in
  $ L_k^n $ containing this cone are exactly $ \sigma_{\{0,\dots,k-2,i\}} $ for
  $ i=k-1,\dots,n $. Consequently, the balancing condition that we have to
  check is just
    \[ u_{k-1} + \cdots + u_n = 0
         \quad \in \RR^n / (\RR \, u_0 + \cdots + \RR \, u_{k-2}), \]
  which is obvious since $ u_0 + \cdots + u_n = 0 $. So the $ L_k^n $ are
  tropical fans of dimension $k$ in $ \RR^n $. We can think of them as
  the ``standard $k$-dimensional linear subspace in $ \RR^n $''.
\end {example}

\begin {definition}[Irreducible fans] \label {def-irred}
  Let $X$ be a tropical fan in $V$. We say that $X$ is \df {irreducible} if
  there is no tropical fan $Y$ of the same dimension in $V$ with $ |Y|
  \subsetneq |X| $.
\end {definition}

\begin {remark} \label {rem-irred}
  Note that the condition of irreducibility remains unchanged under
  refinements, i.e.\ it is well-defined on equivalence classes of tropical
  fans.
\end {remark}

\begin {example} \label {ex-irred}
  Tropical lines, i.e.\ the tropical fans $ L_1^n $ of example \ref {ex-fan}
  \ref {ex-fan-c} (see also example \ref {ex-marked}) are irreducible: any
  1-dimensional fan $Y$ in $ \RR^n $ with $ |Y| \subsetneq |L_1^n| $ would have
  to be obtained by simply removing some of the edges of $ L_1^n $ (and
  possibly changing the weights at the remaining ones), and it is clear that
  doing so would spoil the balancing condition at the origin since the
  remaining edge vectors are linearly independent. (In fact, all linear spaces
  $ L_k^n $ of example \ref {ex-fan} \ref {ex-fan-c} are irreducible, but we
  will not prove this here as we will not need this result.)

  In the same way it follows of course that $ \RR^1 $ (and in fact any tropical
  fan that is just a straight line through the origin in some $ \RR^n $) is
  irreducible.
\end {example}

\begin {remark} \label {rem-irred-2}
  Although there is a notion of irreducible tropical fans it should be noted
  that there is nothing like a unique decomposition into irreducible
  components. The easiest example for this is probably the following
  1-dimensional tropical fan in $ \RR^2 $ (where the vectors in the picture
  denote the markings as in definition \ref {def-marked}).
  \begin {center} \begin{picture}(0,0)%
\includegraphics{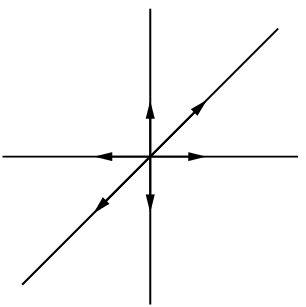}%
\end{picture}%
\setlength{\unitlength}{4144sp}%
\begingroup\makeatletter\ifx\SetFigFont\undefined%
\gdef\SetFigFont#1#2#3#4#5{%
  \reset@font\fontsize{#1}{#2pt}%
  \fontfamily{#3}\fontseries{#4}\fontshape{#5}%
  \selectfont}%
\fi\endgroup%
\begin{picture}(1374,1374)(439,-748)
\end{picture}%
 \end {center}
  By example \ref {ex-irred} it can be realized as a union of irreducible
  tropical fans either as a union of three straight lines through the origin,
  or as a union of the tropical line $ L_1^2 $ of example \ref {ex-fan} \ref
  {ex-fan-c} and the same line mirrored at the origin.
\end {remark}

\begin {proposition} \label {prod-irred}
  The product of two irreducible tropical fans (see example \ref {ex-tropfan}
  \ref {ex-tropfan-c}) is irreducible.
\end {proposition}

\begin {proof}
  Let $X$ and $X'$ be irreducible tropical fans of dimensions $n$ and $n'$
  in $ V = \Lambda \otimes \RR $ and $ V' = \Lambda' \otimes \RR $,
  respectively, and let $Y$ be a tropical fan of dimension $ n+n' $ in $ V
  \times V' $ with $ |Y| \subsetneq |X \times X'| $. By passing from $Y$ to the
  refinement $ Y \cap (X \times X') $ (see example \ref {ex-refine} \ref
  {ex-refine-b} and \ref {ex-refine-d}) we may assume that $ Y \subset X \times
  X' $. For any $ \sigma \times \sigma' \in (X \times X')^{(n+n')} $ we claim
  that $ A(\sigma \times \sigma') := C_{Y,X \times X'}^{-1} (\sigma \times
  \sigma') \cap Y^{(n+n')} $ is either empty or consists of cones that cover $
  \sigma \times \sigma' $ and all have the same weight. In fact, if this was
  not the case then there would have to be a cone $ \tau \in Y^{(n+n'-1)} $
  with $ C_{Y,X \times X'}(\tau) = \sigma \times \sigma' $ so that to the two
  sides of $ \tau $ in $ \sigma \times \sigma' $ there is either only one cone
  of $ A(\sigma \times \sigma') $ or two cones with different weights --- and
  in both cases the balancing condition for $Y$ would be violated at $ \tau $.

  Let us now define a purely $ (n+n') $-dimensional fan $Z$ in $ V \times V' $
  by taking all cones $ \sigma \times \sigma' \in (X \times X')^{(n+n')} $ for
  which $ A(\sigma \times \sigma') \neq \emptyset $, together with all faces of
  these cones. Associating to such a cone $ \sigma \times \sigma'
  \in Z^{(n+n')} $ the weight of the cones in $ A(\sigma \times \sigma') $ we
  make $Z$ into a tropical fan of dimension $ n+n' $ --- in fact it is just a
  tropical fan of which $Y$ is a refinement.

  As $ \dim Z = n+n' $ and $ |Z| = |Y| \subsetneq |X \times X'| $ there must be
  cones $ \sigma_1 \times \sigma_1' \in Z $ and $ \sigma_0 \times \sigma_0' \in
  (X \times X') \backslash Z $. We distinguish two cases:
  \begin {itemize}
  \item $ \sigma_0 \times \sigma_1' \in (X \times X') \backslash Z $: Then we
    construct a purely $n$-dimensional fan $ X_0 $ in $V$ by taking all cones $
    \sigma \in X^{(n)} $ with $ \sigma \times \sigma_1' \in Z $, together with
    all faces of these cones. Setting $ \omega_{X_0} (\sigma) := \omega_Z
    (\sigma \times \sigma_1') $ the fan $ X_0 $ becomes an $n$-dimensional
    tropical fan, with the balancing condition inherited from $Z$ (in fact, the
    balancing condition around a cone $ \tau \in X_0^{(n-1)} $ follows from the
    one around $ \tau \times \sigma_1' $ in $Z$). As $ |X_0| $ is neither empty
    (since $ \sigma_1 \in X_0 $) nor all of $ |X| $ (since $ \sigma_0 \notin
    X_0 $) this is a contradiction to $X$ being irreducible.
  \item $ \sigma_0 \times \sigma_1' \in Z $: This case follows in the same way
    by considering the purely $n'$-dimensional fan $ Y_0 $ in $V'$ given by all
    cones $ \sigma' \in Y^{(n')} $ with $ \sigma_0 \times \sigma' \in Z $,
    together with all faces of these cones --- leading to a contradiction to
    $X'$ being irreducible.
  \end {itemize}
\end {proof}

\begin {lemma} \label {sub-irred}
  Let $X$ and $Y$ be tropical fans of dimension $n$ in $V$. Assume that $ |Y|
  \subset |X| $, and that $X$ is irreducible. Then $ Y \cong \lambda \cdot X $
  for some $ \lambda \in \QQ_{>0} $ (see example \ref {ex-tropfan} \ref
  {ex-tropfan-b}).
\end {lemma}

\begin {proof}
  As $X$ is irreducible we have $ |Y|=|X| $. Replacing $X$ and $Y$ by
  the refinements $ X \cap Y $ and $ Y \cap X $ respectively (see example \ref
  {ex-refine} \ref {ex-refine-b}) we may assume that $X$ and $Y$ consist of the
  same cones (with possibly different weights). Let $ \lambda :=
  \min_{\sigma \in X^{(n)}} \omega_Y(\sigma) / \omega_X (\sigma) > 0 $, and let
  $ \alpha \in \ZZ_{>0} $ with $ \alpha \lambda \in \ZZ $. Consider the new
  weight function $ \omega (\sigma) = \alpha (\omega_Y (\sigma) - \lambda \,
  \omega_X (\sigma)) $ for $ \sigma \in X^{(n)} $. By construction it takes
  values in $ \ZZ_{\ge 0} $, with value 0 occurring at least once. Construct a
  new weighted fan $Z$ from this weight function by taking all cones $ \sigma
  \in X^{(n)} $ with $ \omega(\sigma)>0 $, together with all faces of these
  cones. Then $ (Z,\omega) $ is a tropical fan as the balancing condition
  is linear in the weights. It does not cover $|X|$ since at least one maximal
  cone has been deleted from $X$ by construction. Hence it must be empty as $X$
  was assumed to be irreducible. This means that $ \omega_Y (\sigma) - \lambda
  \, \omega_X (\sigma) = 0 $ for all $ \sigma \in X^{(n)} $.
\end {proof}

\begin {definition}[Morphisms of fans] \label {def-morphism}
  Let $X$ be a fan in $ V = \Lambda \otimes \RR $, and let $Y$ be a fan in $ V'
  = \Lambda' \otimes \RR $. A \df {morphism} $ f: X \to Y $ is a map $f$ from $ |X| \subset V $ to $ |Y| \subset V' $ such that $f$ is $ \ZZ
  $-linear, i.e.\ induced by a linear map from $ \Lambda $ to $ \Lambda' $. By abuse of
  notation we will often denote the corresponding linear maps from $ \Lambda $
  to $ \Lambda' $ and from $V$ to $V'$ by the same letter $f$ (note however
  that these latter maps are not determined uniquely by $ f: X \to Y $ if $ |X|
  $ does not span $V$). A morphism of weighted fans is a morphism of fans
  (i.e.\ there are no conditions on the weights).
\end {definition}

\begin {remark} \label {rem-morphism}
  It is obvious from the definition that morphisms from $ (X,\omega_X)
  $ to $ (Y,\omega_Y) $ are in one-to-one correspondence with morphisms 
 from any refinement of $  (X,\omega_X) $ to any refinement of
  $ (Y,\omega_Y) $ simply by using the same map $f:|X| \rightarrow |Y|$.
\end {remark}

\begin {construction} \label {constr-image}
  Let $X$ be a purely $n$-dimensional fan in $ V = \Lambda \otimes \RR $, and
  let $Y$ be any fan in $ \Lambda' \otimes \RR $. For any morphism $ f: X \to Y
  $ we will construct an \df {image fan} $ f(X) $ in $V'$ (which is empty or of pure dimension $n$)
  as follows. Consider the collection of cones
    \[ Z := \{ f(\sigma) ;\; \mbox {$ \sigma \in X $ contained in a maximal
               cone of $X$ on which $f$ is injective} \} \]
  in $V'$. Note that $Z$ is in general not a fan in $V'$: it satisfies
  condition \ref {def-fan-a} of definition \ref {def-fan}, but in general not
  \ref {def-fan-b} (since e.g.\ the images of some maximal cones might
  overlap in a region of dimension $n$). To make it into one choose linear
  forms $ g_1,\dots,g_N \in \Lambda'^\vee $ such that each cone $ f(\sigma) \in Z $
  can be written as
    \[ f(\sigma) = \{ x \in V' ;\; \mbox {$ g_i(x) = 0 $, $ g_j(x) \ge 0 $,
         $ g_k(x) \le 0 $ for all $ i \in I_\sigma, j \in J_\sigma,
         k \in K_\sigma $} \} \tag {$*$} \]
  for suitable index sets $ I_\sigma, J_\sigma, K_\sigma \in \{1,\dots,N\} $.
  Now we replace $X$ by the fan $ \tilde X = X \cap H_{g_1 \circ f} \cap \cdots
  \cap H_{g_N \circ f} $ (see example \ref {ex-fan} \ref {ex-fan-b}). This fan
  satisfies $ |\tilde X|=|X| $, and by definition each cone of $ \tilde X $ is
  of the form
    \[ \{ x \in V ;\; \mbox {$ x \in \sigma $, $ g_i(f(x)) = 0 $, $ g_j(f(x))
         \ge 0 $, $ g_k(f(x)) \le 0 $ for all $ i \in I, j \in J, k \in K $}
         \} \]
  for some $ \sigma \in X $ and a partition $ I \cup J \cup K = \{1,\dots,N\}
  $. By $(*)$ the image of such a cone under $f$ is of the form
  \begin {align*}
    & \{ x \in V' ;\; \mbox {$ x \in f(\sigma) $, $ g_i(x) = 0 $, $ g_j(x) \ge
         0 $, $ g_k(x) \le 0 $ for all $ i \in I, j \in J, k \in K $} \} \\
    & \quad = \{ x \in V' ;\; \mbox {$ g_i(x) = 0 $, $ g_j(x) \ge 0
         $, $ g_k(x) \le 0 $ for all $ i \in I', j \in J', k \in K' $} \}
  \end {align*}
  for some partition $ I' \cup J' \cup K' = \{1,\dots,N\} $, i.e.\ it is a cone
  of the fan $ H_{g_1} \cap \cdots \cap H_{g_N} $. Hence the set of cones
    \[ \tilde Z := \{ f(\sigma) ;\; \mbox {$ \sigma \in \tilde X $ contained
               in a maximal cone of $ \tilde X $ on which $f$ is injective} \}
       \]
  consists of cones in the fan $ H_{g_1} \cap \cdots \cap H_{g_N} $. As these
  cones now satisfy condition \ref {def-fan-b} of definition \ref {def-fan} it
  follows that $ \tilde Z $ is a fan in $V'$. It is obvious that it is of pure
  dimension $n$.

  If moreover $X$ is a weighted fan then $ \tilde X $ will be a weighted fan as
  well (see example \ref {ex-refine} \ref {ex-refine-b}). In this case we make
  $ \tilde Z $ into a weighted fan by setting
    \[ \omega_{\tilde Z} (\sigma') :=
         \sum_{\sigma \in \tilde X^{(n)}: f(\sigma) = \sigma'}
           \omega_{\tilde X}(\sigma) \cdot
             |\Lambda'_{\sigma'}/f(\Lambda_\sigma)| \]
  for all $ \sigma' \in \tilde Z^{(n)} $ (see definition \ref {def-cone}; note
  that $ f(\Lambda_\sigma) $ is a sublattice of $ \Lambda'_{\sigma'} $ of full
  rank).

  We denote the fan $ (\tilde Z,\omega_{\tilde Z}) $ by $ f(X) $ and call it the
  \df {image fan} of $f$.
  It is clear from the construction that the equivalence class of $f(X)$
  remains unchanged if we replace the $ g_1,\dots,g_N $
  by any larger set of linear forms (this would just lead to a refinement), and
  hence also if we replace them by any other set of linear forms satisfying
  $(*)$. Hence the equivalence class of $ f(X) $ does
  not depend on any choices we made. 
  It is obvious that the equivalence class of $ f(X)
  $ depends in fact only on the equivalence class of $X$.

  By abuse of notation we will usually drop the tilde from the above notation
  and summarize the construction above as follows: given a weighted fan $X$ of
  dimension $n$, an arbitrary fan $Y$, and a morphism $ f:X \to Y $, we may
  assume after passing to an equivalent fan for $X$ that
    \[ f(X) = \{ f(\sigma) ;\; \mbox {$ \sigma \in X $ contained in a maximal
               cone of $X$ on which $f$ is injective} \} \]
  is a fan. We consider this to be a weighted fan of dimension $n$ by setting
    \[ \omega_{f(X)} (\sigma') :=
         \sum_{\sigma \in X^{(n)}: f(\sigma) = \sigma'}
           \omega_X (\sigma) \cdot |\Lambda'_{\sigma'}/f(\Lambda_\sigma)| \]
  for $ \sigma' \in f(X)^{(n)} $. This weighted fan is well-defined up to
  equivalence.
\end {construction}

\begin {proposition} \label {prop-image}
  Let $X$ be an $n$-dimensional tropical fan in $ V = \Lambda \otimes \RR $,
  $Y$ an arbitrary fan in $ V' = \Lambda' \otimes \RR $, and $ f: X \to Y $ a
  morphism. Then $ f(X) $ is an $n$-dimensional tropical fan as well (if it is
  not empty).
\end {proposition}

\begin {proof}
  Since we have already seen in construction \ref {constr-image} that $ f(X) $
  is a weighted fan of dimension $n$ we just have to check the balancing
  condition of definition \ref {def-tropfan}. As in this construction we may
  assume that $ f(X) $ just consists of the cones $ f(\sigma) $ for all $
  \sigma \in X $ contained in a maximal cone on which $f$ is injective. Let $
  \tau' \in f(X)^{(n-1)} $, and let $ \tau \in X^{(n-1)} $ with $ f(\tau) =
  \tau' $. Applying $f$ to the balancing condition for $X$ at $ \tau $ we get
    \[ \sum_{\sigma > \tau} \omega_X(\sigma) \cdot f(u_{\sigma/\tau}) = 0
         \quad \in V' / V'_{\tau'} \]
  for all $ \tau \in X^{(n-1)} $. Now let $ \sigma' \in f(X) $ be a cone with $
  \sigma' > \tau' $, and $ \sigma \in X $ with $ \sigma > \tau $ such that $
  f(\sigma) = \sigma' $. Note that the primitive normal vector $
  u_{\sigma'/\tau'} $ is related to the (possibly non-primitive) vector $
  f(u_{\sigma/\tau}) $ by
    \[ f(u_{\sigma/\tau})
         = |\Lambda'_{\sigma'} / (\Lambda'_{\tau'} + \ZZ\,f(u_{\sigma/\tau})) |
           \cdot u_{\sigma'/\tau'} \]
  if $f$ is injective on $ \sigma $, and $ f(u_{\sigma/\tau}) = 0 $ otherwise.
  Inserting this into the above balancing condition, and using the exact
  sequence
    \[ 0 \;\longrightarrow\; \Lambda'_{\tau'}/f(\Lambda_\tau)
         \;\longrightarrow\; \Lambda'_{\sigma'}/f(\Lambda_\sigma)
         \;\longrightarrow\; \Lambda'_{\sigma'}/(\Lambda'_{\tau'} +
             \ZZ \, f(u_{\sigma/\tau}))
         \;\longrightarrow\; 0 \]
  we conclude that
    \[ \sum_{\sigma>\tau} \omega_X(\sigma)
         \cdot |\Lambda'_{\sigma'}/f(\Lambda_\sigma)|
         \cdot u_{\sigma'/\tau'} = 0 \quad \in V' / V'_{\tau'}, \]
  where the sum is understood to be taken over only those $ \sigma>\tau $ on
  which $f$ is injective.

  Let us now sum these equations up for all $ \tau $ with $ f(\tau) = \tau' $.
  The above sum then simply becomes a sum over all $ \sigma $ with $
  f(\sigma)>\tau' $ (note that each such $ \sigma $ occurs in the sum exactly
  once since $f$ is injective on $ \sigma $ so that $ \sigma $ cannot have two
  distinct codimension-1 faces that both map to $ \tau' $). Splitting this sum
  up according to the cone $ f(\sigma) $ we get
    \[ \sum_{\sigma'>\tau'} \;\;\; \sum_{\sigma: f(\sigma) = \sigma'}
         \omega_X(\sigma)
         \cdot |\Lambda'_{\sigma'}/f(\Lambda_\sigma)|
         \cdot u_{\sigma'/\tau'} = 0 \quad \in V' / V'_{\tau'}. \]
  But using the definition of the weights of $ f(X) $ of construction \ref
  {constr-image} this is now simply the balancing condition
    \[ \sum_{\sigma'>\tau'} \omega_{f(X)}(\sigma')
         \cdot u_{\sigma'/\tau'} = 0 \quad \in V' / V'_{\tau'} \]
  for $ f(X) $.
\end {proof}

\begin {corollary} \label {cor-image}
  Let $X$ and $Y$ be tropical fans of the same dimension $n$ in $ V = \Lambda
  \otimes \RR $ and $ V' = \Lambda' \otimes \RR $, respectively, and let $ f: X
  \to Y $ be a morphism. Assume that $Y$ is irreducible. Then there is a fan
  $ Y_0 $ in $V'$ of smaller dimension with $ |Y_0| \subset |Y| $ such that
  \begin {enumerate}
  \item \label {cor-image-a}
    each point $ Q \in |Y| \backslash |Y_0| $ lies in the interior of a
    cone $ \sigma_Q' \in Y $ of dimension $n$;
  \item \label {cor-image-b}
    each point $ P \in f^{-1} (|Y| \backslash |Y_0|) $ lies in the interior
    of a cone $ \sigma_P \in X $ of dimension $n$;
  \item \label {cor-image-c}
    for $ Q \in |Y| \backslash |Y_0| $ the sum
      \[ \sum_{P \in |X|: f(P)=Q} \mult_P f \]
    does not depend on $Q$, where the multiplicity $ \mult_P f $ of $f$ at $P$
    is defined to be
      \[ \mult_P f := \frac {\omega_X(\sigma_P)}{\omega_Y(\sigma'_Q)}
           \cdot |\Lambda'_{\sigma'_Q}/f(\Lambda_{\sigma_P})| \]
  \end {enumerate}
\end {corollary}

\begin {proof}
  Consider the tropical fan $ f(X) $ in $V'$ (see construction \ref
  {constr-image} and proposition \ref {prop-image}). If $ f(X) = \emptyset $
  (i.e.\ if there is no maximal cone of $X$ on which $f$ is injective) the
  statement of the corollary is trivial. Otherwise $ f(X) $ has dimension $n$
  and satisfies $ |f(X)| \subset |Y| $, so as $Y$ is irreducible it follows by
  lemma \ref {sub-irred} that $ f(X) \cong \lambda \cdot Y $ for some $ \lambda
  \in \QQ_{>0} $. After passing to equivalent fans we may assume that $ f(X) $
  and $Y$ consist of the same cones, and that these are exactly the cones of
  the form $ f(\sigma) $ for $ \sigma \in X $ contained in a maximal cone on
  which $f$ is injective. Now let $ Y_0 $ be the fan consisting of all cones of
  $Y$ dimension less than $n$. Then \ref {cor-image-a} and \ref {cor-image-b}
  hold by construction. Moreover, each $ Q \in |Y| \backslash |Y_0| $ lies in
  the interior of a unique $n$-dimensional cone $ \sigma' $, and there is a
  1:1 correspondence between points $ P \in f^{-1}(Q) $ and $n$-dimensional
  cones $ \sigma $ in $X$ with $ f(\sigma) = \sigma' $. So we conclude
  that
    \[ \sum_{P: f(P)=Q} \mult_P f
       = \sum_{\sigma: f(\sigma) = \sigma'}
           \frac {\omega_X (\sigma)}{\omega_Y(\sigma')}
           \cdot |\Lambda'_{\sigma'}/f(\Lambda_\sigma)|
       = \frac {\omega_{f(X)} (\sigma')}{\omega_Y(\sigma')}
       = \lambda \]
  does not depend on $Q$.
\end {proof}

\begin {corollary} \label {cor-det}
  In the situation and with the notation of corollary \ref {cor-image} assume
  moreover that $X$ and $Y$ are marked fans as in definition \ref {def-marked},
  and that their structure of tropical fans has been induced by these data as
  in construction \ref {constr-marked} \ref {constr-marked-a}. Let $
  \sigma_1,\dots,\sigma_n \in X^{(1)} $ be the 1-dimensional faces of $
  \sigma_P $, and let $ \sigma'_1,\dots,\sigma'_n \in Y^{(1)} $ be the
  1-dimensional faces of $ \sigma'_Q = \sigma'_{f(P)} $. Then $ \mult_P f $ is
  equal to the absolute value of the determinant of the matrix for the linear
  map $ f|_{V_{\sigma_P}}: V_{\sigma_P} \to V'_{\sigma'_Q} $ in the bases $ \{
  v_{\sigma_1},\dots,v_{\sigma_n} \} $ and $ \{ v_{\sigma_1'},\dots,
  v_{\sigma_n'} \} $.
\end {corollary}

\begin {proof}
  It is well-known that $ |\Lambda'_{\sigma'_Q}/f(\Lambda_{\sigma_P})| $ is the
  determinant of the linear map $ f|_{V_{\sigma_P}}: V_{\sigma_P} \to
  V'_{\sigma'_Q} $ with respect to lattice bases of $ \Lambda_{\sigma_P} $ and
  $ \Lambda'_{\sigma'_Q} $. The statement of the corollary now follows since
  the base change to $ \{ v_{\sigma_1},\dots,v_{\sigma_n} \} $ and $ \{
  v_{\sigma_1'},\dots,v_{\sigma_n'} \} $ clearly leads to factors in the
  determinant of $ \omega_X(\sigma_P) $ and $ 1/\omega_Y(\sigma'_Q) $,
  respectively.
\end {proof}

  \section{Tropical $\calM_{0,n}$, Grassmannians, and the space of trees}
  \label{sec-m0n}

In this chapter we want to show that the moduli space $\calM_{0,n,\trop}$ of
rational $n$-marked abstract tropical curves is a tropical fan in the sense of
definition \ref{def-tropfan}, and that the forgetful map is a morphism of fans.

Let us start by recalling the relevant definitions from \cite {GM}.

\begin {definition}[Graphs] \label {def-graph} 
  \begin {enumerate}
  \item \label {def-graph-a}
    Let $ I_1,\dots,I_n \subset \RR $ be a finite set of closed, bounded or
    half-bounded real intervals. We pick some (not necessarily distinct)
    boundary points $ P_1,\dots,P_k,Q_1,\dots,Q_k \in I_1 \dunion \cdots
    \dunion I_n $ of these intervals. The topological space $ \Gamma $ obtained
    by identifying $ P_i $ with $ Q_i $ for all $ i=1,\dots,k $ in $ I_1
    \dunion \cdots \dunion I_n $ is called a \df {graph}. As usual, the \df
    {genus} of $ \Gamma $ is simply its first Betti number $ \dim H_1
    (\Gamma,\RR) $.
  \item \label {def-graph-b}
    For a graph $ \Gamma $ the boundary points of the intervals $ I_1,\dots,
    I_n $ are called the \df {flags}, their image points in $ \Gamma $ the \df
    {vertices} of $ \Gamma $. If $F$ is such a flag then its image vertex in $
    \Gamma $ will be denoted $ \partial F $. For a vertex $V$ the number of
    flags $F$ with $ \partial F = V $ is called the \df {valence} of $V$ and
    denoted $ \val V $. We denote by $ \Gamma^0 $ and $ \Gamma' $ the sets of
    vertices and flags of $ \Gamma $, respectively.
  \item \label {def-graph-c}
    The open intervals $ I_1^\circ,\dots,I_n^\circ $ are naturally open subsets
    of $ \Gamma $; they are called the \df {edges} of $ \Gamma $. An edge will
    be called bounded (resp.\ unbounded) if its corresponding open interval is.
    We denote by $ \Gamma^1 $ (resp.\ $ \Gamma^1_0 $ and $ \Gamma^1_\infty $)
    the set of edges (resp.\ bounded and unbounded edges) of $ \Gamma $. Every
    flag $ F \in \Gamma' $ belongs to exactly one edge that we will denote by $
    [F] \in \Gamma^1 $. The unbounded edges will also be called the \df {ends}
    of $ \Gamma $.
  \end {enumerate}
\end {definition}

\begin {definition}[Abstract tropical curves] \label {def-tropcurve}
  A (rational, abstract) tropical curve is a connected graph $ \Gamma $ of
  genus 0 all of whose vertices have valence at least 3. An \df {$n$-marked
  tropical curve} is a tuple $ (\Gamma,x_1,\dots,x_n) $ where $ \Gamma $ is a
  tropical curve and $ x_1,\dots,x_n \in \Gamma^1_\infty $ are distinct
  unbounded edges of $ \Gamma $. Two such marked tropical curves $ (\Gamma,
  x_1,\dots,x_n) $ and $ (\tilde \Gamma, \tilde x_1,\dots,\tilde x_n) $ are
  called isomorphic (and will from now on be identified) if there is a
  homeomorphism $ \Gamma \to \tilde \Gamma $ mapping $ x_i $ to $ \tilde x_i $
  for all $i$ and such that every edge of $ \Gamma $ is mapped bijectively onto
  an edge of $ \tilde \Gamma $ by an affine map of slope $ \pm 1 $, i.e.\ by a
  map of the form $ t \mapsto a \pm t $ for some $ a \in \RR $. The space of
  all $n$-marked tropical curves (modulo isomorphisms) with precisely $n$
  unbounded edges will be denoted $ \calM_{0,n,\trop} $. (It can be thought of as a
  tropical analogue of the moduli space $ \bar M_{0,n} $ of $n$-pointed stable
  rational curves.)
\end {definition}

% 
% 
% A (rational) \df {abstract tropical curve} is a connected rational graph
% $\Gamma$ (that is, a tree) whose vertices have valence at least $3$. Unbounded
% edges (also called ends or leaves) are allowed. The bounded edges are equipped
% with a positive length. An $n$-marked abstract tropical curve is a tuple
% $(\Gamma,x_1,\ldots,x_n)$ where $\Gamma$ is an abstract tropical curve and
% $x_1,\ldots,x_n$ are distinct unbounded edges. For a more detailed definition
% of abstract tropical curve, see \cite{GM} definition 2.2.

% \begin{definition} \label{def-m0n}
%   The space $ \calM_{0,n,\trop} $ is defined to be the space of all $n$-marked
%   abstract tropical curves (modulo isomorphisms) with exactly $n$ leaves.
% \end{definition}

\begin {remark} \label {can-par}
  The isomorphism condition of definition \ref {def-tropcurve} means that every
  edge of a marked tropical curve has a parametrization as an interval in $ \RR
  $ that is unique up to translations and sign. In particular, every bounded
  edge $E$ of a tropical curve has an intrinsic \df {length} that we will
  denote by $ l(E) \in \RR_{>0} $.

  One way to fix this translation and sign ambiguity is to pick a flag $F$ of
  the edge $E$: there is then a unique choice of parametrization such that the
  corresponding closed interval is $ [0,l(E)] $ (or $ [0,\infty) $ for
  unbounded edges), with the chosen flag $F$ being the zero point of this
  interval. We will call this the \df {canonical parametrization} of $E$ with
  respect to the flag $F$.
\end {remark}
\begin {definition}[Combinatorial types]
  The \df {combinatorial type} of a marked tropical curve $ (\Gamma,x_1,\dots,
  x_n) $ is defined to be the homeomorphism class of $ \Gamma $ relative $
  x_1,\dots,x_n $ (i.e.\ the data of $ (\Gamma,x_1,\dots,x_n) $ modulo
  homeomorphisms of $ \Gamma $ that map each $ x_i $ to itself).
\end {definition}

% The \df {(combinatorial) type} of a marked abstract tropical curve
% $(\Gamma,x_1,\ldots,x_n)$ is the data left when dropping the information about
% the lengths of the bounded edges. 
Lemma 2.10 of \cite{GM} says that there are
only finitely many combinatorial types of curves in $\calM_{0,n}$. The subset of curves of type $\alpha$
in $ \calM_{0,n,\trop} $ is the interior of a cone in $\RR^k$ (where $k$ is the
number of bounded edges) given by the inequalities that all lengths are
positive --- i.e.\ it is the positive orthant of $\RR^k$. Example 2.13 of
\cite{GM} describes how these cones are glued locally in $ \calM_{0,n,\trop} $.

Note that this construction is exactly the same as the construction of the
\df {space of (phylogenetic) trees} $\mathbb{T}_n$ in \cite{BHV}, page 9--13.

Fix $ n \ge 3 $ and consider the space $\RR^{\binom{n}{2}}$ indexed by the
set $\mathcal{T}$ of all subsets $T \subset [n] := \{1,\dots,n\} $ with $
|T|=2 $. In order to embed $\calM_{0,n,\trop}$ into a quotient of
$\RR^{\binom{n}{2}}$, consider the map
\begin{eqnarray*}
  \dist_n: \calM_{0,n,\trop} & \longrightarrow  & \RR^{\binom{n}{2}} \\
           (\Gamma,x_1,\ldots,x_n) & \longmapsto
             & (\dist_\Gamma(x_i,x_j))_{\{i,j\}\in\mathcal{T}}
\end{eqnarray*}
where $\dist_\Gamma(x_i,x_j)$ denotes the distance between the unbounded edges
(or leaves) $x_i$ and $x_j$, that is, the sum of the lengths of all edges on
the (unique) path leading from $x_i$ to $x_j$. Furthermore, define a linear map
$\Phi_n$ by
\begin{eqnarray*}
 \Phi_n: \RR^n & \longrightarrow & \RR^{\binom{n}{2}}\\
          a & \longmapsto & (a_i+a_j)_{{\{i,j\}\in\mathcal{T}}}.
\end{eqnarray*}
Denote by $Q_n$ the $(\binom{n}{2}-n)$-dimensional quotient vector space
$\RR^{\binom{n}{2}}/\im(\Phi_n)$, and by $q_n:\RR^{\binom{n}{2}}\rightarrow
Q_n$ the canonical projection.

\begin{theorem}[\cite{SS}, theorem 4.2] \label{SpeyerSturmfelsTheorem}
  The map $ \varphi_n := q_n \circ \dist_n:\calM_{0,n,\trop}\rightarrow Q_n$
  is an embedding, and the image $\varphi_n(\calM_{0,n,\trop})\subset Q_n$ is a
  simplicial fan of pure dimension $n-3$. The interior of its $k$-dimensional
  cells corresponds to combinatorial types of graphs with $n$ marked leaves and
  exactly $k$ bounded edges.
\end{theorem}

\begin{construction} \label {constr-vi}
  For each subset $I\subset [n]$ of cardinality $1<|I|<n-1$, define $v_I$ to be
  the image under $\varphi_n$ of a tree with one bounded edge of length one,
  the marked ends with labels in $I$ on one side of the bounded edge and the
  marked ends with labels in $[n]\setminus I$ on the other. Note that $ v_I =
  v_{[n] \backslash I} $ by construction. As an example, the following picture
  shows the 5-marked rational curve corresponding to the vector $ v_{\{1,2,5\}}
  = v_{\{3,4\}} \in Q_5 $:

  \begin {center} \begin{picture}(0,0)%
\includegraphics{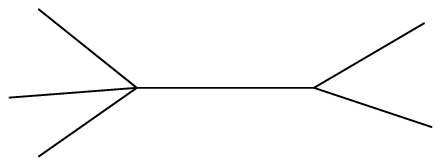}%
\end{picture}%
\setlength{\unitlength}{4144sp}%
\begingroup\makeatletter\ifx\SetFigFontNFSS\undefined%
\gdef\SetFigFontNFSS#1#2#3#4#5{%
  \reset@font\fontsize{#1}{#2pt}%
  \fontfamily{#3}\fontseries{#4}\fontshape{#5}%
  \selectfont}%
\fi\endgroup%
\begin{picture}(2064,916)(698,-965)
\put(864,-196){\makebox(0,0)[rb]{\smash{{\SetFigFontNFSS{10}{12.0}{\familydefault}{\mddefault}{\updefault}{\color[rgb]{0,0,0}$x_1$}%
}}}}
\put(713,-608){\makebox(0,0)[rb]{\smash{{\SetFigFontNFSS{10}{12.0}{\familydefault}{\mddefault}{\updefault}{\color[rgb]{0,0,0}$x_2$}%
}}}}
\put(864,-901){\makebox(0,0)[rb]{\smash{{\SetFigFontNFSS{10}{12.0}{\familydefault}{\mddefault}{\updefault}{\color[rgb]{0,0,0}$x_5$}%
}}}}
\put(2738,-257){\makebox(0,0)[lb]{\smash{{\SetFigFontNFSS{10}{12.0}{\familydefault}{\mddefault}{\updefault}{\color[rgb]{0,0,0}$x_3$}%
}}}}
\put(2747,-736){\makebox(0,0)[lb]{\smash{{\SetFigFontNFSS{10}{12.0}{\familydefault}{\mddefault}{\updefault}{\color[rgb]{0,0,0}$x_4$}%
}}}}
\end{picture}%
 \end {center}

  By theorem \ref{SpeyerSturmfelsTheorem}, the $v_I$ generate the edges of the
  simplicial fan $\varphi_n(\calM_{0,n,\trop})$. Let $\Lambda_n:=\langle v_I
  \rangle_\ZZ\subset Q_n$ be the lattice in $Q_n$ generated by the vectors
  $v_I$.
\end{construction}

\begin{theorem} \label{thm-m0nfan}
  The marked fan $(\varphi_n(\calM_{0,n,\trop}),\{v_I\})$ is a tropical fan
  in $Q_n$ with lattice $\Lambda_n$. In other words, using the embedding $
  \varphi_n $ the space $ \calM_{0,n,\trop} $ can be thought of as a tropical
  fan of dimension $ n-3 $.
\end{theorem}

\begin{proof}
  By lemma \ref{lem-marked} we have to prove that
    \[ \sum_{\sigma>\tau} v_{\sigma/\tau} = 0\in  Q_n/V_\tau \]
  for all $\tau\in\varphi_n(\calM_{0,n,\trop})^{(n-4)}$ and $v_{\sigma/\tau}$
  as in construction \ref {constr-marked} \ref {constr-marked-b}.

  In order to prove this fix a cell $\tau\in\varphi_n(\calM_{0,n,\trop})^{
  (n-4)}$. As $\tau$ is a cell of codimension one, the curves in the interior
  of $ \tau $ have exactly one vertex $V$ of valence 4 and all other vertices
  of valence 3. Denote by $ I_s \subset [n] $ for $ s=1,\dots,4 $ the set of
  marked ends lying behind the $s$-th edge adjacent to $V$. By construction
  there are exactly three combinatorial types $ \sigma_1,\sigma_2,\sigma_3>\tau
  $ in $ \calM_{0,n,\trop} $, obtained from $\tau$ by replacing the vertex $V$
  with an edge as in the following picture (where we have only drawn the
  relevant part of the graphs and the labels denote the marked ends lying
  behind the corresponding edges):
  
  \begin {center} \begin{picture}(0,0)%
\includegraphics{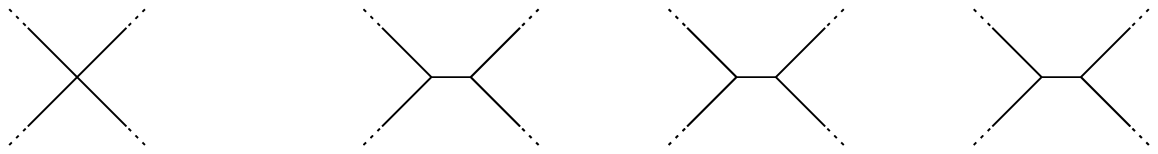}%
\end{picture}%
\setlength{\unitlength}{4144sp}%
\begingroup\makeatletter\ifx\SetFigFontNFSS\undefined%
\gdef\SetFigFontNFSS#1#2#3#4#5{%
  \reset@font\fontsize{#1}{#2pt}%
  \fontfamily{#3}\fontseries{#4}\fontshape{#5}%
  \selectfont}%
\fi\endgroup%
\begin{picture}(5340,1066)(976,-1160)
\put(991,-241){\makebox(0,0)[rb]{\smash{{\SetFigFontNFSS{10}{12.0}{\familydefault}{\mddefault}{\updefault}{\color[rgb]{0,0,0}$I_1$}%
}}}}
\put(991,-871){\makebox(0,0)[rb]{\smash{{\SetFigFontNFSS{10}{12.0}{\familydefault}{\mddefault}{\updefault}{\color[rgb]{0,0,0}$I_2$}%
}}}}
\put(1711,-871){\makebox(0,0)[lb]{\smash{{\SetFigFontNFSS{10}{12.0}{\familydefault}{\mddefault}{\updefault}{\color[rgb]{0,0,0}$I_3$}%
}}}}
\put(1711,-241){\makebox(0,0)[lb]{\smash{{\SetFigFontNFSS{10}{12.0}{\familydefault}{\mddefault}{\updefault}{\color[rgb]{0,0,0}$I_4$}%
}}}}
\put(2611,-241){\makebox(0,0)[rb]{\smash{{\SetFigFontNFSS{10}{12.0}{\familydefault}{\mddefault}{\updefault}{\color[rgb]{0,0,0}$I_1$}%
}}}}
\put(2611,-871){\makebox(0,0)[rb]{\smash{{\SetFigFontNFSS{10}{12.0}{\familydefault}{\mddefault}{\updefault}{\color[rgb]{0,0,0}$I_2$}%
}}}}
\put(3511,-241){\makebox(0,0)[lb]{\smash{{\SetFigFontNFSS{10}{12.0}{\familydefault}{\mddefault}{\updefault}{\color[rgb]{0,0,0}$I_4$}%
}}}}
\put(3511,-871){\makebox(0,0)[lb]{\smash{{\SetFigFontNFSS{10}{12.0}{\familydefault}{\mddefault}{\updefault}{\color[rgb]{0,0,0}$I_3$}%
}}}}
\put(3061,-1096){\makebox(0,0)[b]{\smash{{\SetFigFontNFSS{10}{12.0}{\familydefault}{\mddefault}{\updefault}{\color[rgb]{0,0,0}$\sigma_1$}%
}}}}
\put(4006,-241){\makebox(0,0)[rb]{\smash{{\SetFigFontNFSS{10}{12.0}{\familydefault}{\mddefault}{\updefault}{\color[rgb]{0,0,0}$I_1$}%
}}}}
\put(4906,-241){\makebox(0,0)[lb]{\smash{{\SetFigFontNFSS{10}{12.0}{\familydefault}{\mddefault}{\updefault}{\color[rgb]{0,0,0}$I_4$}%
}}}}
\put(5401,-241){\makebox(0,0)[rb]{\smash{{\SetFigFontNFSS{10}{12.0}{\familydefault}{\mddefault}{\updefault}{\color[rgb]{0,0,0}$I_1$}%
}}}}
\put(5401,-871){\makebox(0,0)[rb]{\smash{{\SetFigFontNFSS{10}{12.0}{\familydefault}{\mddefault}{\updefault}{\color[rgb]{0,0,0}$I_4$}%
}}}}
\put(4006,-871){\makebox(0,0)[rb]{\smash{{\SetFigFontNFSS{10}{12.0}{\familydefault}{\mddefault}{\updefault}{\color[rgb]{0,0,0}$I_3$}%
}}}}
\put(4906,-871){\makebox(0,0)[lb]{\smash{{\SetFigFontNFSS{10}{12.0}{\familydefault}{\mddefault}{\updefault}{\color[rgb]{0,0,0}$I_2$}%
}}}}
\put(6301,-871){\makebox(0,0)[lb]{\smash{{\SetFigFontNFSS{10}{12.0}{\familydefault}{\mddefault}{\updefault}{\color[rgb]{0,0,0}$I_2$}%
}}}}
\put(6301,-241){\makebox(0,0)[lb]{\smash{{\SetFigFontNFSS{10}{12.0}{\familydefault}{\mddefault}{\updefault}{\color[rgb]{0,0,0}$I_3$}%
}}}}
\put(1366,-443){\makebox(0,0)[b]{\smash{{\SetFigFontNFSS{10}{12.0}{\familydefault}{\mddefault}{\updefault}{\color[rgb]{0,0,0}$V$}%
}}}}
\put(1351,-1096){\makebox(0,0)[b]{\smash{{\SetFigFontNFSS{10}{12.0}{\familydefault}{\mddefault}{\updefault}{\color[rgb]{0,0,0}$\tau$}%
}}}}
\put(4456,-1096){\makebox(0,0)[b]{\smash{{\SetFigFontNFSS{10}{12.0}{\familydefault}{\mddefault}{\updefault}{\color[rgb]{0,0,0}$\sigma_2$}%
}}}}
\put(5851,-1096){\makebox(0,0)[b]{\smash{{\SetFigFontNFSS{10}{12.0}{\familydefault}{\mddefault}{\updefault}{\color[rgb]{0,0,0}$\sigma_3$}%
}}}}
\end{picture}%
 \end {center}

  Note that $ v_{I_1 \cup I_2} $ is an edge of $ \sigma_1 $ which is not an
  edge of $ \tau $; hence we have $ v_{\sigma_1/\tau} = v_{I_1\cup I_2} $, and
  similarly $ v_{\sigma_2/\tau} = v_{I_1 \cup I_3} $ and $ v_{\sigma_3/\tau} =
  v_{I_1 \cup I_4} $.

  Now let $ C \in \tau $ be the point corresponding to the curve obtained from
  the combinatorial type $\tau$ by setting the lengths of each bounded edge $E$
  to 1 if $E$ is adjacent to $V$, and 0 otherwise:

  \begin {center} \begin{picture}(0,0)%
\includegraphics{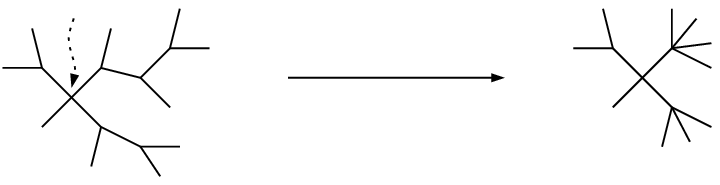}%
\end{picture}%
\setlength{\unitlength}{4144sp}%
\begingroup\makeatletter\ifx\SetFigFontNFSS\undefined%
\gdef\SetFigFontNFSS#1#2#3#4#5{%
  \reset@font\fontsize{#1}{#2pt}%
  \fontfamily{#3}\fontseries{#4}\fontshape{#5}%
  \selectfont}%
\fi\endgroup%
\begin{picture}(3264,1295)(979,-1331)
\put(1336,-159){\makebox(0,0)[b]{\smash{{\SetFigFontNFSS{10}{12.0}{\familydefault}{\mddefault}{\updefault}{\color[rgb]{0,0,0}$V$}%
}}}}
\put(3916,-1141){\makebox(0,0)[b]{\smash{{\SetFigFontNFSS{10}{12.0}{\familydefault}{\mddefault}{\updefault}{\color[rgb]{0,0,0}$C$}%
}}}}
\put(1441,-1096){\makebox(0,0)[b]{\smash{{\SetFigFontNFSS{10}{12.0}{\familydefault}{\mddefault}{\updefault}{\color[rgb]{0,0,0}a curve in the}%
}}}}
\put(1441,-1276){\makebox(0,0)[b]{\smash{{\SetFigFontNFSS{10}{12.0}{\familydefault}{\mddefault}{\updefault}{\color[rgb]{0,0,0}interior of $\tau$}%
}}}}
\end{picture}%
 \end {center}

Note that this operation changes the combinatorial type in general.
  Furthermore, let $ a \in \RR^n $ be the vector with $ a_i = 1 $ if the
  marked end $ x_i $ is adjacent to $V$, and $ a_i=0 $ otherwise. We claim
  that
    \[ v_{\sigma_1/\tau} + v_{\sigma_2/\tau} + v_{\sigma_3/\tau}
         = \Phi_n(a) + \dist_n(C) \in \RR^{\binom n2}, \]
  from which the required balancing condition then follows after passing to
  the quotient by $ \im \Phi_n + V_\tau $. We check this equality
  coordinate-wise for all $T=\{i,j\}\in\mathcal{T}$. We may assume that $i$
  and $j$ do not lie in the same $ I_s $ since otherwise the $T$-coordinate of
  every term in the equation is zero. Then the $T$-coordinate on the left hand
  side is 2 since $ x_i $ and $ x_j $ are on different sides of the newly
  inserted edge for exactly two of the types $ \sigma_1,\sigma_2,\sigma_3 $.
  On the other hand, if we denote by $ 0\leq c \leq 2 $ the number how many of
  the two ends $ x_i $ and $ x_j $ are adjacent to $V$, then the
  $T$-coordinates of $ \Phi_n(a) $ and $ \dist_n(C) $ are $c$ and $ 2-c $
  respectively; so the proposition follows.
\end{proof}

\begin {remark}
  A slightly different proof of this balancing condition has recently been
  found independently by Mikhalkin in \cite {M2} section 2.
\end {remark}

\begin{remark}
  The space of trees appears in \cite{SS} as a quotient of the tropical
  Grassmannian. Let us review how the tropical Grassmannian is defined.
  Tropical varieties can be defined as images of usual algebraic varieties over
  the field $K$ of Puiseux series. An element of $K$ is a Puiseux series $p(t)=
  a_1 t^{q_1}+a_2t^{q_2}+ a_3t^{q_3} +\ldots$, where $a_i \in \CC$ and $q_1
  <q_2 <q_3< \ldots$ are rational numbers which share a common denominator.
  Define the \df {valuation} of $p(t)$ to be $ \val(p(t))=-q_1$, if $p(t)\neq
  0$, and $-\infty$ otherwise. Let $V(I)$ be a variety in $(K^*)^n$. Its
  \df {tropicalization} $\Trop(I)$ is defined to be the closure of the set
    \[ \{\big(\val(y_1),\ldots,\val(y_n)\big)\;|\;
         (y_1,\ldots,y_n)\in V(I)\}\subset \RR^n. \]
  The tropicalization of an ideal can also be computed using initial ideals: it
  is equal to the closure of the set of $w\in \QQ^n$ such that the initial
  ideal $\ini_w I$ contains no monomial (see theorem 2.1 of \cite{SS}). In this
  way, it can be considered as a subfan of the Gr\"obner fan of $I$ and
  inherits its fan structure (this is explained for example in \cite{BJSST}).
  Theorem 2.5.1 of \cite{S} shows that $\Trop(I)$ is a tropical fan.

  In \cite{SS}, the \df {tropical Grassmannian} $\mathcal{G}_{2,n}$ is defined
  to be the tropicalization of the ideal $I_{2,n}$ of the usual Grassmannian in
  its Pl\"ucker embedding. The ideal $I_{2,n}$ is an ideal in the polynomial ring with
  $\binom{n}{2}$ variables. For example, we have
    \[ I_{2,4}=\langle
         p_{12}p_{34}-p_{13}p_{24}+p_{14}p_{23}\rangle\subset
         K[p_{12},p_{13},p_{14},p_{23},p_{24},p_{34}]. \]
  (We label the variables with subsets of size $2$ of $\{1,\ldots,n\}$.)
  The tropical Grassmannian $\mathcal{G}_{2,n}:= \Trop(I_{2,n})$ is a tropical
  fan in $\RR^{\binom{n}{2}}$.
\end{remark}

The lineality space of $\mathcal{G}_{2,n}$ (that is, the intersection of all
cones of $\mathcal{G}_{2,n}$) is equal to the image of $\Phi_n$ (see \cite{SS},
page 6). As $\im(\Phi_n)$ is contained in every cone of $\mathcal{G}_{2,n}$, we
can reduce $\mathcal{G}_{2,n}$ modulo $\im(\Phi_n)$ and define the
\df {reduced tropical Grassmannian}
  \[ \mathcal{G}_{2,n}':=\mathcal{G}_{2,n}/\im(\Phi_n)\subset
       \RR^{\binom{n}{2}}/\im(\Phi_n). \]
Then $\mathcal{G}_{2,n}'$ is a tropical fan, too. Theorem 3.4 of \cite{SS}
states that $\mathcal{G}_{2,n}'$ is equal to the space of phylogenetic trees
$\mathbb{T}_n$.

\begin{example} \label {ex-m04}
  For $ n=4 $ the space $ \calM_{0,4,\trop} $ is embedded by $ \varphi_4 $ in
  $ Q_4 = \RR^{\binom 42} / \im \Phi_n \cong \RR^6 / \RR^4 \cong \RR^2 $, and
  the lattice $ \Lambda_4 \subset Q_4 $ is spanned by the three vectors $
  v_{\{1,2\}} = v_{\{3,4\}} $, $ v_{\{1,3\}} = v_{\{2,4\}} $, and $
  v_{\{1,4\}} = v_{\{2,3\}} $. The space $ \calM_{0,4,\trop} $ is of dimension
  1 and has three maximal cones, each spanned by one of the three vectors
  above. By the balancing condition of theorem \ref {thm-m0nfan} the sum of
  these three vectors is 0 in $ Q_4 $, and hence any two of them form a basis
  of $ \Lambda_4 $. With (the negative of) such a basis $ \varphi_4(
  \calM_{0,4,\trop}) $ simply becomes the standard tropical line $ L_1^2
  \subset \RR^2 $ as in example \ref {ex-fan} \ref {ex-fan-c}.
\end{example}

For the rest of this chapter we will now consider the forgetful maps between
the moduli spaces of abstract curves and show that they are morphisms of fans.
For simplicity we will only consider the map that forgets the last marked end.

\begin{definition}[Forgetful maps] \label{def-forget}
   Let $n \geq 4$ be an integer. We have a \df {forgetful map} $\ft$ from
   $\calM_{0,n,\trop}$ to $\calM_{0, n - 1,\trop}$ which assigns to an
   $n$-marked curve $(C,x_1,\dots ,x_n)$ the (stabilization of the)
   $(n-1)$-marked curve $(C,x_1 ,\dots ,x_{n-1})$ (see \cite{GM}, definition
   4.1). By stabilization we mean that we have to straighten $2$-valent vertices that possibly emerge when removing the marked end $x_n$. 
\end{definition}

\begin{proposition}\label{prop-forgetmorphism}
  With the tropical fan structure of theorem \ref {thm-m0nfan} the forgetful
  map $\ft: \calM_{0,n,\trop}\rightarrow \calM_{0,n-1,\trop}$ is a morphism
  of fans in the sense of definition \ref{def-morphism}.
\end{proposition}

\begin{proof}
  Let $\pr : \RR^{\binom{n}{2}} \rightarrow \RR^{\binom{n - 1}{2}}$
  denote the projection to those coordinates $T\in\mathcal{T}$ with $n
  \not\in T$. As $\pr(\im(\Phi_n))=\im(\Phi_{n-1})$, the map $\pr$ induces a
  linear map $\widetilde{\pr} : Q_n\rightarrow Q_{n-1}$.

  We claim that $\widetilde{\pr}_{|\varphi_n(M_{0,n,\trop})}$ is the map
  induced by $\ft$, hence we have to show the commutativity of the diagram
    \[ \begin{xy}
       \xymatrix{ \calM_{0,n,\trop} \ar[r]^-{\ft} & \calM_{0,n-1,\trop} \\
         Q_n \ar[r]^-{\widetilde{\pr}} \ar@{<-}[u]^{\varphi_{n}} &
         Q_{n-1}. \ar@{<-}[u]_{\varphi_{n-1}}}
       \end{xy} \]
  So let $C = (\Gamma, x_1, \ldots ,x_n) \in \calM_{0,n,\trop}$ be an abstract
  $n$-marked rational tropical curve. If $x_n$ is not adjacent to exactly one
  bounded and one unbounded edge, then forgetting $x_n$ does not change any of
  the distances $\dist_C(x_i,x_j)$ for $i,j \neq n$; so in this case we are
  done.

  Now assume that $x_n$ is adjacent to one unbounded edge $x_k$ and one bounded
  edge $E$ of length $l(E)$. Then the distances between marked points are given
  by
    \[ \dist_{\ft(C)}(x_i,x_j)=
       \begin{cases}
         \dist_C(x_i,x_j)-l(E) & \text{ if } j=k\\
         \dist_C(x_i,x_j)      & \text{ if } j\neq k,
       \end{cases} \]
  hence $\widetilde{\pr} (\varphi_n (C)) = \varphi_{n - 1} (\ft (C)) + \Phi_{
  n-1}(l (E) \cdot e_k)$, where $e_k$ denotes the $k$-th standard unit vector
  in $\RR^{n-1}$. We conclude that $\widetilde{\pr} \circ \varphi_n =
  \varphi_{n - 1} \circ \ft$. Hence the diagram is commutative.

  It remains to check that $\widetilde{\pr}(\Lambda_n)\subset \Lambda_{n-1}$.
  It suffices to check this on the generators $v_I$; and we can assume $n\in
  I$ since $v_I=v_{[n]\setminus I} $. We get
    \[ \widetilde{\pr}(v_I) = \begin{cases}
         v_{I\setminus\{n\}} & \text{, if } |I|\geq 3 \\
         0 & \text {, if } |I|=2.
       \end{cases} \]
\end{proof}

  \section{Tropical curves in $ \RR^r $} \label {sec-smap}

In this chapter we will consider tropical analogues of the algebro-geometric
moduli spaces of (rational) stable maps as e.g.\ already introduced in \cite
{GM}. Similarly to the classical case the points of these spaces correspond to
marked abstract tropical curves together with a suitable map to some $ \RR^r $.
We would like to make these spaces into tropical fans by considering the
underlying abstract tropical curves and using the fan structure of $
\calM_{0,n,\trop} $ developed in chapter \ref {sec-m0n}. In order for this to
work it turns out however that we have to modify the construction of \cite {GM}
slightly: in the underlying abstract tropical curves we have to label \textsl
{all} ends and not just the contracted ones corresponding to the marked points.

\begin{definition}[Tropical $ \calM_{0,n,\trop}^{\lab}(\RR^r,\Delta) $]
    \label{def-m0nprd}
  A \df {(parametrized) labeled $n$-marked tropical curve in} $ \RR^r $ is a
  tuple $ (\Gamma,x_1,\dots,x_N,h) $ for some $ N \ge n $, where $
  (\Gamma,x_1,\dots,x_N) $ is an abstract $N$-marked tropical curve with exactly $N$ ends and $ h:
  \Gamma \to \RR^r $ is a continuous map satisfying:
  \begin {enumerate}
  \item On each edge of $ \Gamma $ the map $h$ is of the form $ h(t) = a + t
    \cdot v $ for some $ a \in \RR^r $ and $ v \in \ZZ^r $. The integral
    vector $v$ occurring in this equation if we pick for $E$ the canonical
    parametrization starting at $V\in \partial E$ (see \ref{can-par})
    will be denoted $ v(E,V) $ and called the \df {direction} of $E$ (at
    $V$). If $E$ is an unbounded edge and $V$ its only boundary point we will
    write for simplicity $ v(E) $ instead of $ v(E,V) $.
  \item For every vertex $V$ of $ \Gamma $ we have the \df {balancing
    condition}
      \[ \sum_{E| V \in \partial E} v(E,V) = 0. \]
  \item $ v(x_i)=0 $ for $ i=1,\dots,n $ (i.e.\ each of the first $n$ ends is
    contracted by $h$), whereas $ v(x_i) \neq 0 $ for $ i>n $ (i.e.\ the
    remaining $ N-n $ ends are ``non-contracted ends'').
  \end {enumerate}
  Two labeled $n$-marked tropical curves $ (\Gamma,x_1,\dots,x_N,h) $ and $
  (\tilde \Gamma, \tilde x_1,\dots,\tilde x_N,\tilde h)$ in $ \RR^r $ are
  called isomorphic (and will from now on be identified) if there is an
  isomorphism $\varphi: (\Gamma,x_1,\dots,x_N) \to (\tilde \Gamma,\tilde x_1,
  \dots,\tilde x_N) $ of the underlying abstract curves such that $ \tilde h
  \circ \varphi = h $.

  The \df {degree} of a labeled $n$-marked tropical curve as above is defined
  to be the $(N-n)$-tuple $ \Delta = (v(x_{n+1}),\dots,v(x_N)) \in (\ZZ^r
  \backslash \{0\})^{N-n}$ of directions of its non-contracted ends. 
 Its \df
  {combinatorial type} is given by the data of the combinatorial type of the
  underlying abstract marked tropical curve $ (\Gamma,x_1,\dots,x_N) $
  together with the directions of all its (bounded as well as unbounded)
  edges. For the rest of this work the number $N$ will always be related to $n$
  and $ \Delta $ by $ N=n+\#\Delta $ and thus denote the total number of
  (contracted or non-contracted) ends of an $n$-marked curve in $ \RR^r $ of
  degree $ \Delta $.
  
  The space of all labeled $n$-marked tropical curves of a given degree
  $\Delta$ in $\RR^r$ will be denoted $ \calM_{0,n,\trop}^{\lab}(\mathbb{R}^r,
  \Delta) $. For the special choice
    \[ \Delta = (-e_0,\dots,-e_0\;\;,\dots,\;\;-e_r,\dots,-e_r) \]
  with $ e_0 := -e_1-\cdots-e_r $ and where each $ e_i $ occurs
  exactly $d$ times we will also denote this space by $ \calM_{0,n,\trop}^{
  \lab}(\mathbb{R}^r,d) $ and say that these curves have degree $d$.
\end{definition}

\begin{definition}[Evaluation map] \label{def-ev}
  For $ i=1,\dots,n $ the map
  \begin{eqnarray*}
    \ev_i:
      \calM_{0,n,\trop}^{\text{lab}}(\mathbb{R}^r,\Delta)&\rightarrow&\RR^r\\
      (\Gamma,x_1,\dots x_N,h) & \longmapsto & h(x_i)
  \end{eqnarray*}
  is called the {\em $i$-th evaluation map} (note that this is well-defined
  for the contracted ends since for them $ h(x_i) $ is a point in $ \RR^r $).
\end{definition}

\begin{construction}[Tropical $ \calM_{0,n,\trop}(\RR^r,\Delta) $]
    \label{constr-m0nrd}
  Let $ N \ge n \ge 0 $, and let $ \Delta =(v_{n+1},\dots,v_N) \in (\ZZ^r
  \backslash \{0\})^{N-n} $. It is obvious from the definitions that the
  subgroup $G$ of the symmetric group $ \SM_{N-n} $ consisting of all
  permutations $ \sigma $ of $ \{ n+1,\dots,N \} $ such that $ v_{\sigma(i)}
  = v_i $ for all $ i=n+1,\dots,N $ acts on the space $ \calM_{0,n,\trop}^{\lab}
  (\mathbb{R}^r,\Delta) $ by relabeling the non-contracted ends. We denote the
  quotient space $ \calM_{0,n,\trop}^{\lab}(\mathbb{R}^r,\Delta) / G $ by $
  \calM_{0,n,\trop}(\mathbb{R}^r,\Delta) $ and call this the space of \df
  {(unlabeled) $n$-marked tropical curves} in $ \RR^r $ of degree $ \Delta $
  in $ \RR^r $; its elements can obviously be thought of as $n$-marked
  tropical curves in $ \RR^r $ for which we have only specified how many of
  its ends have a given direction, but where no labeling of these
  (non-contracted) ends is given. Consequently, when considering unlabeled
  curves we can (and usually will) think of $ \Delta $ as a multiset
  $ \{ v_{n+1},\dots,v_N \} $ instead of as a vector $ (v_{n+1},\dots,v_N)
  $.

  Note that this definition of $ \calM_{0,n,\trop}(\mathbb{R}^r,\Delta) $
  agrees with the one given in \cite {GM}, and that for curves of degree $d$
  the group $G$ above is precisely $ (\SM_d)^{r+1} $.
\end{construction}

\begin{remark}
  Proposition 2.1 of \cite{NS} tells us that there are only finitely many
  combinatorial types in any given moduli space $ \calM_{0,n,\trop}(\RR^r,
  \Delta) $ (and thus also in $ \calM_{0,n,\trop}^{\lab}(\RR^r,\Delta)$).
  Proposition 2.11 of \cite{GM} says that the subset of curves of a fixed
  combinatorial type is a cone in a real vector space, and example 2.13
  shows how these cones are glued locally. (These results are only stated there
  for plane tropical curves, but it is obvious that they hold in the same way
  if we replace $\RR^2$ by $\RR^r$.) Hence the moduli spaces of labeled or
  unlabeled $n$-marked tropical curves in $ \RR^r $ are polyhedral complexes.
\end{remark}

\begin{remark}
  We will now show how one can use the tropical fan structure of
  $\calM_{0,n,\trop}$ to make the moduli spaces $ \calM_{0,n,\trop}^{\lab}(
  \mathbb{R}^r,d)$ of labeled curves in $ \RR^r $ into tropical fans as well.
  The idea is the same as in the classical algebro-geometric case: the
  non-contracted ends of the tropical curves in $ \RR^r $ can be thought of as
  their intersection points with the $ r+1 $ coordinate hyperplanes of $ \PP^r
  $; so passing from unlabeled to labeled curves corresponds to labeling these
  intersection points (which is one possible strategy to construct the moduli
  spaces of stable maps in algebraic geometry). The final moduli spaces of
  stable maps then arise by taking the quotient by the group of possible
  permutations of the labels. Note that this makes the spaces of stable maps
  into stacks instead of varieties (as the group action is in general not
  free); and similarly our moduli spaces of unlabeled tropical curves in $
  \RR^r $ can only be thought of as ``tropical stacks'' instead of tropical
  varieties. To avoid this notion of tropical stacks (which has not been
  developed yet) we will always work with labeled tropical curves in this
  paper. Of course this is not a problem for enumerative purposes since we can
  always count labeled curves first and then divide the result by $ |G| $ in
  the end.
\end{remark}

Now let us make the relation between the spaces $ \calM_{0,n,\trop}^{\lab}(
\mathbb{R}^r,\Delta)$ and $\calM_{0,N,\trop}$ (with $ N = n+\#\Delta $)
precise. For this we consider the forgetful map
\begin {align*}
  \psi: \quad
        \calM_{0,n,\trop}^{\lab}(\mathbb{R}^r,\Delta) &\to
          \calM_{0,n+\#\Delta,\trop} \\
          (\Gamma,x_1,\dots,x_N,h) &\mapsto
          (\Gamma,x_1,\dots,x_N)
\end {align*}
that forgets the map to $ \RR^r $ but keeps all (contracted and non-contracted)
unbounded edges. For this construction (and the rest of this paper) we will assume
for simplicity that $ N \ge 3 $, so that this map is well-defined.

\begin{lemma} \label {lem-combtypebij}
  The forgetful map $ \psi $ induces a bijection of combinatorial types of the
  two moduli spaces $ \calM_{0,n,\trop}^{\lab}(\mathbb{R}^r,\Delta) $ and $
  \calM_{0,N,\trop} $ with $ N=n+\#\Delta $.
\end{lemma}

\begin{proof}
  It is clear that $ \psi $ induces a well-defined map between combinatorial
  types of labeled and abstract tropical curves. By definition the only
  additional data in the combinatorial type of a labeled curve compared to the
  underlying abstract curve is the directions of the edges. As the directions
  of the ends are fixed by $ \Delta $ it therefore suffices to show that for a
  given (combinatorial type of a) graph $ \Gamma $ and fixed directions of the
  ends there is a unique choice of directions for the bounded edges compatible
  with the balancing condition of definition \ref {def-m0nprd}.

  We will prove this by induction on $N$. The statement is clear for $ N=3 $
  since there are no bounded edges in this case. In order to prove the
  induction step we show first that there is a vertex $V_1$ of valence $l$
  with $l-1$ ends adjacent. Assume there is no such vertex. Then
    \[ N\leq \sum_V (\val V-2)=\sum_V \val V-2 \cdot k, \]
  where $k$ denotes the number of vertices. As
    \[ k=N-2-\sum_V (\val V-3)=N-2-\sum_V \val V+3k \]
  we have
    \[ N\leq \sum_V \val V- \sum_V \val V +N-2 \]
  which is a contradiction.

  So assume that the (contracted or non-contracted) ends $x_{i_1},\ldots,
  x_{i_{l-1}} $ are adjacent at a vertex $V_1$ of valence $l$. Let $E$ be the
  only bounded edge which is adjacent to $V_1$. Then the directions of
  $x_{i_1},\ldots,x_{i_{l-1}} $ are fixed (by $ \Delta $ for the non-contracted
  ends, and to be 0 for the contracted ones), so there is a unique choice for
  the direction of $E$ compatible with the balancing condition at $V_1$. Now
  remove the unbounded edges $x_{i_1},\ldots,x_{i_{l-1}}$ from $\Gamma$ and
  make $E$ into an unbounded edge with this new fixed direction. This new
  graph has $ N-(l-1)+1 $ ends, which is less than $N$ as $l$ is at least 3.
  Therefore we can apply induction and conclude that the directions of all
  bounded edges are determined.
\end{proof}

For the rest of this paper we will assume for simplicity that $ n>0 $, i.e.\
that there is at least one contracted end.

\begin {proposition} \label{prop-m0nrdfan}
  With notations as above, the map
  \begin {align*}
    \calM_{0,n,\trop}^{\lab}(\mathbb{R}^r,\Delta)
      &\to Q_N \times \RR^r \\
    C &\mapsto (\varphi_N(\psi(C)),\ev_1(C))
  \end {align*}
  is an embedding whose image is the tropical fan $ \varphi_N (\calM_{0,N,
  \trop}) \times \RR^r $. So $ \calM_{0,n,\trop}^{\lab}(\mathbb{R}^r,\Delta) $
  can be thought of as a tropical fan of dimension $ r+N-3=r+n+\#\Delta-3 $,
  namely as the fan $ \calM_{0,n+\#\Delta,\trop} \times \RR^r $.
\end {proposition}

\begin{proof}
  It is clear that the given map is a continuous map of polyhedral complexes.
  By lemma \ref {lem-combtypebij} it suffices to check injectivity and compute
  its image for a fixed combinatorial type $\alpha$. As in proposition 2.11 of \cite{GM} the cell of $ \calM_{0,n,\trop}^{\lab}(
  \mathbb{R}^r,\Delta) $ corresponding to curves of type $ \alpha $ is given by
  $ \RR_{>0}^k \times \RR^r $, where the first positive coordinates are the
  lengths of the bounded edges, and the second $r$ coordinates are the position
  of a root vertex (that we choose to be the vertex $V$ adjacent to $ x_1 $ here). 
This is true because we can recover the map $h$ if we know the position $h(V)$ and the lengths of each bounded edge (in addition to the information of $\alpha$).
The map of the
  proposition obviously sends this cell bijectively to the product of the
  corresponding cell of $ \varphi_N (\calM_{0,N,\trop})$ and $ \RR^r $ (since
  the position of the root vertex $h(V)$ is unrestricted).
\end{proof}

\begin{proposition} \label{prop-evmorphism}
  With the tropical fan structure of proposition \ref {prop-m0nrdfan} the
  evaluation maps $ \ev_i: \calM_{0,n,\trop}^{\lab}(\mathbb{R}^r,\Delta)
  \rightarrow \RR^r$ are morphisms of fans (in the sense of definition \ref
  {def-morphism}).
\end{proposition}

\begin{proof}
  As usual let $ N=n+\#\Delta $, and identify $\calM_{0,n,\trop}^{\text{lab}}(
  \mathbb{R}^r,\Delta)$ with the space $\calM_{0,N,\trop}\times \RR^r$ as in
  proposition \ref{prop-m0nrdfan}. For all $1\leq i\leq n$ consider the linear
  map
  \begin{eqnarray*}
     \ev_i':  \RR^{\binom{N}{2}} \times \RR^r & \longrightarrow & \RR^r\\
     (a_{\{1,2\}},\dots, a_{\{N-1,N\}},b)  & \longmapsto &
     b + \frac{1}{2}\sum_{\substack {k=2 \\ k\neq i}}^{N}
       \left(a_{\{1,k\}}-a_{\{i,k\}}\right)v_k,
  \end{eqnarray*}
  where $ v_k $ is the direction of the $k$-th end, i.e.\ $ v_1 = \cdots = v_n
  = 0 $ and $ (v_{n+1},\dots,v_N) = \Delta $. As $\ev_i'(\im(\Phi_N)\times
  \{0\})=\{0\}$, the map $\ev_i'$ induces a linear map $\widetilde{\ev_i'} :
  Q_N\times\RR^r\rightarrow \RR^r$.

  We claim that $\widetilde{\ev_i'}|_{\calM_{0,N,\trop} \times \RR^r}$ is the
  map induced by $\ev_i$, hence we have to show the commutativity of the
  diagram
    \[ \begin{xy}
       \xymatrix{ \calM_{0,N,\trop}\times\RR^r \ar[r]^{\qquad \ev_{i}}
         & \RR^r \\
       Q_N \times \RR^r \ar@{<-}[u]^{\varphi_{N}\times\id}
         \ar[ur]^{\widetilde{\ev_{i}'}}}
       \end{xy} \]
  As the case $i=1$ is trivial, and using the additivity of the function
  $\sum_{k=1}^{N} v_k\left(a_{1k}-a_{ik}\right)$, we may assume that there
  exists only one edge $E$, and that this edge separates $x_1$ and $x_i$. Let
  $V$ be the vertex adjacent to $ x_1 $, and let $ T_1 \subset [N] $ be the
  set of ends lying on the connected component of $ \Gamma\setminus
  \{E\}$ containing $x_1$ (so that $ 1 \in T $ and $ i \notin T $). Since by
  the balancing condition of definition \ref{def-m0nprd} the equation $
  v(E,V)=\sum_{k\notin T_1}v_k=-\sum_{k\in T_1}v_k$ holds, we get
  \begin{eqnarray*}
    \left(\widetilde{\ev_i'}\circ(\varphi_N\times{id})\right)(\Gamma,x_1,
      \dots, x_N,h(x_1))
      & = & h(x_1)+\frac{1}{2}\left[-\sum_{k\in T_1}l(E)v_k
          +\sum_{k\not\in T_1}l(E)v_k\right] \\
      & = & h(x_1)+l(E) \, v(E,V) \\
      & = & h(x_i)
  \end{eqnarray*}
  as required (where $ l(E) $ denotes the length of $E$). It remains to check
  that $\widetilde{\ev_i'}(\Lambda_N\times\ZZ^r)\subset \ZZ^r$. It suffices
  to check that $\widetilde{\ev_i'}(v_I,b)\in \ZZ^r$ for the generators $v_I$
  of $\Lambda_N$ (see construction \ref {constr-vi}) and all $b\in\ZZ^r$.
  Assuming $1\in I$, we get
    \[ \widetilde{\ev_i'}(v_I,b) = \begin{cases}
         b - \sum_{k\in I} v_k & \text{, if } i\not\in I\\
         b & \text {, if } i\in I,
       \end{cases} \]
  which finishes the proof as $v_k\in\ZZ^r$.
\end{proof}

\begin{example}
  As an example of the calculation in the above proof we consider the space $
  \calM_{0,2,\trop}^{\lab} (\RR^2,1) = \calM_{0,5,\trop} \times \RR^2 $, so
  that the curves have $N=5$ unbounded edges with directions $ v_1=v_2=0 $, $
  v_3=(1,1) $, $ v_4 = (-1,0) $, $ v_5 = (0,-1) $. We consider curves of the
  combinatorial type drawn in the following picture:

  \begin {center} \begin{picture}(0,0)%
\includegraphics{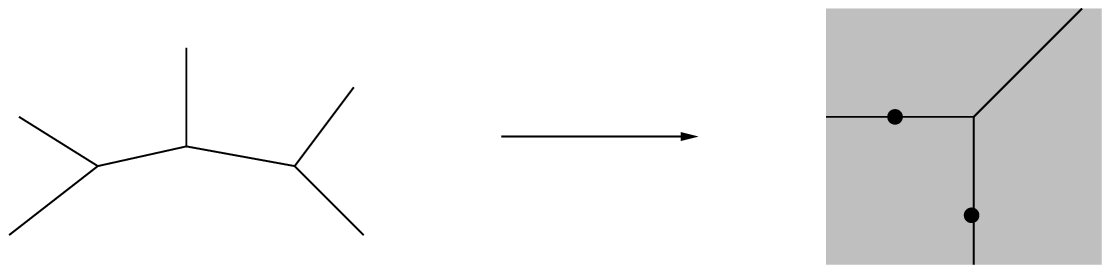}%
\end{picture}%
\setlength{\unitlength}{4144sp}%
\begingroup\makeatletter\ifx\SetFigFontNFSS\undefined%
\gdef\SetFigFontNFSS#1#2#3#4#5{%
  \reset@font\fontsize{#1}{#2pt}%
  \fontfamily{#3}\fontseries{#4}\fontshape{#5}%
  \selectfont}%
\fi\endgroup%
\begin{picture}(5111,1606)(3365,-2915)
\put(6121,-1996){\makebox(0,0)[b]{\smash{{\SetFigFontNFSS{10}{12.0}{\familydefault}{\mddefault}{\updefault}{\color[rgb]{0,0,0}$h$}%
}}}}
\put(7471,-1906){\makebox(0,0)[b]{\smash{{\SetFigFontNFSS{10}{12.0}{\familydefault}{\mddefault}{\updefault}{\color[rgb]{0,0,0}$h(x_1)$}%
}}}}
\put(7876,-2491){\makebox(0,0)[lb]{\smash{{\SetFigFontNFSS{10}{12.0}{\familydefault}{\mddefault}{\updefault}{\color[rgb]{0,0,0}$h(x_2)$}%
}}}}
\put(7651,-2122){\makebox(0,0)[b]{\smash{{\SetFigFontNFSS{10}{12.0}{\familydefault}{\mddefault}{\updefault}{\color[rgb]{0,0,0}$a$}%
}}}}
\put(7858,-2266){\makebox(0,0)[lb]{\smash{{\SetFigFontNFSS{10}{12.0}{\familydefault}{\mddefault}{\updefault}{\color[rgb]{0,0,0}$b$}%
}}}}
\put(3895,-2018){\makebox(0,0)[lb]{\smash{{\SetFigFontNFSS{12}{14.4}{\rmdefault}{\mddefault}{\updefault}{\color[rgb]{0,0,0}   }%
}}}}
\put(3427,-2047){\makebox(0,0)[rb]{\smash{{\SetFigFontNFSS{10}{12.0}{\rmdefault}{\mddefault}{\updefault}{\color[rgb]{0,0,0}$x_1$}%
}}}}
\put(4011,-2129){\makebox(0,0)[b]{\smash{{\SetFigFontNFSS{10}{12.0}{\rmdefault}{\mddefault}{\updefault}{\color[rgb]{0,0,0}$a$}%
}}}}
\put(4545,-2147){\makebox(0,0)[b]{\smash{{\SetFigFontNFSS{10}{12.0}{\rmdefault}{\mddefault}{\updefault}{\color[rgb]{0,0,0}$b$}%
}}}}
\put(5044,-1922){\makebox(0,0)[lb]{\smash{{\SetFigFontNFSS{10}{12.0}{\rmdefault}{\mddefault}{\updefault}{\color[rgb]{0,0,0}$x_2$}%
}}}}
\put(4235,-1603){\makebox(0,0)[b]{\smash{{\SetFigFontNFSS{10}{12.0}{\rmdefault}{\mddefault}{\updefault}{\color[rgb]{0,0,0}$x_3$}%
}}}}
\put(3380,-2607){\makebox(0,0)[rb]{\smash{{\SetFigFontNFSS{10}{12.0}{\rmdefault}{\mddefault}{\updefault}{\color[rgb]{0,0,0}$x_4$}%
}}}}
\put(5093,-2589){\makebox(0,0)[lb]{\smash{{\SetFigFontNFSS{10}{12.0}{\rmdefault}{\mddefault}{\updefault}{\color[rgb]{0,0,0}$x_5$}%
}}}}
\put(8326,-1456){\makebox(0,0)[b]{\smash{{\SetFigFontNFSS{10}{12.0}{\familydefault}{\mddefault}{\updefault}{\color[rgb]{0,0,0}$h(x_3)$}%
}}}}
\put(8461,-1771){\makebox(0,0)[lb]{\smash{{\SetFigFontNFSS{10}{12.0}{\familydefault}{\mddefault}{\updefault}{\color[rgb]{0,0,0}$\RR^2$}%
}}}}
\put(7111,-2041){\makebox(0,0)[rb]{\smash{{\SetFigFontNFSS{10}{12.0}{\familydefault}{\mddefault}{\updefault}{\color[rgb]{0,0,0}$h(x_4)$}%
}}}}
\put(7831,-2851){\makebox(0,0)[b]{\smash{{\SetFigFontNFSS{10}{12.0}{\familydefault}{\mddefault}{\updefault}{\color[rgb]{0,0,0}$h(x_5)$}%
}}}}
\end{picture}%
 \end {center}

  We can read off immediately from the picture that $\ev_2(C)=h(x_1)+a(1,0)+
  b(0,-1)$.

  On the other hand, we compute using the formula of the proof of proposition
  \ref {prop-evmorphism}:
  \begin{align*}
    \left(\widetilde{\ev_2'}\circ(\varphi_N\times{id})\right)
      & (C,x_1,\dots,x_5,h(x_1)) \\
      & = \widetilde{\ev_2'}(\dist_C(x_1,x_2),\dots,\dist_C(x_4,x_5),h(x_1))\\
      & = h(x_1) + \frac{1}{2}\sum_{k=3}^5 \left(\dist_C(x_1,x_k)
          -\dist_C(x_i,x_k)\right)v_k\\
      & = h(x_1) + \frac{1}{2} \left[(a-b)v_3+(0-(a+b))v_4
          + ((a+b)-0)v_5\right]\\
      & = h(x_1) +\frac{1}{2} \left[a(v_3-v_4+v_5) + b(-v_3-v_4+v_5)\right]\\
      & = h(x_1) + a(1,0) + b(0,-1).
  \end{align*}
\end{example}

\begin{remark}
  Just as in chapter \ref{sec-m0n} we can define forgetful maps on the moduli
  space of labeled tropical curves in $ \RR^r $ as well. We can forget certain
  contracted ends, but also we can forget all non-contracted ends and the $
  \RR^r $ factor in $ \calM_{0,n,\trop}^{\lab} (\RR^r,\Delta) = \calM_{0,N,
  \trop} \times \RR^r $, which corresponds to forgetting the map $h$. The same
  argument as in proposition \ref{prop-forgetmorphism} shows that these
  forgetful maps are morphisms as well.
\end{remark}

\begin{remark}
  In proposition \ref {prop-m0nrdfan} we have constructed the tropical fan
  structure on the moduli space $ \calM_{0,n,\trop}^{\lab} (\RR^r,\Delta) $
  using the evaluation at the first (contracted) marked end. If we use a
  different contracted end $ x_i $ with $ i \in \{2,\dots,n\} $ instead, the
  two embeddings in $ Q_N \times \RR^r $ only differ by addition of $ \ev_i(C)
  - \ev_1(C) $ in the $ \RR^r $ factor. As this is an isomorphism by
  proposition \ref {prop-evmorphism} it follows that the tropical fan
  structure defined in proposition \ref {prop-m0nrdfan} is natural in the
  sense that it does not depend on the choice of contracted end.
\end{remark}

  \section{Applications} \label {sec-cor}

In this final chapter we want to apply our results to reprove and generalize
two statements in tropical enumerative geometry that have occurred earlier in
the literature. The first application simply concerns the number of rational
tropical curves in some $ \RR^r $ through given conditions.

\begin{theorem} \label {thm-indep}
  Let $ r \ge 2 $, let $ \Delta $ be a degree of tropical curves in $ \RR^r $
  (i.e.\ a multiset of elements of $ \ZZ^r \backslash \{0\} $ that sum up to $0$), and let $ n>0
  $ be such that $ r+n+\#\Delta-3 = nr $. Then the number of rational tropical
  curves of degree $ \Delta $ in $ \RR^r $ through $n$ points in general
  position (counted with multiplicities) does not depend on the position of
  the points.
\end{theorem}

\begin{proof}
  Let $ \ev:= \ev_1 \times \dots \times \ev_n: \calM_{0,n,\trop}^{\lab}
  (\RR^r,\Delta)\rightarrow \RR^{nr}$. By proposition \ref{prop-m0nrdfan} the
  space $ \calM_{0,n,\trop}^{\lab} (\RR^r,\Delta) $ is a tropical fan, and by
  proposition \ref{prop-evmorphism} the map $\ev$ is a morphism of fans (of the
  same dimension). As $ \RR^{nr} $ is irreducible we know from corollary \ref
  {cor-image} that for general $ Q \in \RR^{nr} $ the number
    \[ \sum_{C \in |\calM_{0,n,\trop}^{\lab}(\RR^r,\Delta)|: \ev(C)=Q}
       \mult_{C} \ev \]
  does not depend on $Q$. So if we define the tropical multiplicity of a curve
  $C$ occurring in this sum to be $ \mult_C \ev $ we conclude that the number
  of labeled curves of degree $ \Delta $ in $ \RR^r $ through the $n$ points
  in $ \RR^r $ specified by $Q$ does not depend on the choice of (general) $Q$.
  The statement for unlabeled curves now follows simply by dividing the
  resulting number by the order $ |G| $ of the symmetry group of the
  non-contracted ends as in construction \ref {constr-m0nrd}.
\end{proof}

\begin{remark} \label {rem-indep}
  In the proof of theorem \ref{thm-indep} we have defined the tropical
  multiplicity of a curve to be the multiplicity of the evaluation map as in
  corollary \ref {cor-image}. By corollary \ref {cor-det} this multiplicity
  can be computed on a fixed cell of $ \calM_{0,n,\trop}^{\lab}(\RR^r,\Delta)
  $ as the absolute value of the determinant of the matrix obtained by
  expressing the evaluation map in terms of the basis on the source given by
  the vectors $ v_I $ (see construction \ref {constr-vi}) that span the given
  cell. But note that the coordinates on the given cell with respect to this
  basis are simply the lengths of the bounded edges. It follows that the
  determinant that we have to consider is precisely the same as that used
  in \cite {GM} chapter 3. In particular, in the case $ r=2 $ it follows by
  proposition 3.8 of \cite {GM} that our multiplicity of the curves agrees with
  Mikhalkin's usual notion of multiplicity as in definitions 4.15 and 4.16 of
  \cite {M0}.
\end{remark}

\begin{remark}
  For simplicity we have formulated theorem \ref{thm-indep} only for the case
  of counting tropical curves through given points. Of course the very same
  proof can be used for counting curves through given affine linear subspaces
  (with rational slopes) if one replaces the evaluation maps $ \ev_i $ by their
  compositions with the quotient maps $ \RR^r \to \RR^r / L_i $, where $ L_i $
  is the linear subspace chosen at the $i$-th contracted end. This setup has
  been considered e.g.\ by Nishinou and Siebert in \cite {NS}.
\end{remark}

As our second application we consider the map occurring in the proof of the
tropical Kontsevich formula (see proposition 4.4 of \cite{GM}).

\begin {theorem} \label {thm-indep2}
  Let $ d \ge 1 $, and let $ n=3d $. Define
    \[ \pi := \ev_1^1 \times \ev_2^2 \times \ev_3 \times \cdots \times \ev_n
       \times \ft_4 : \calM_{0,n,\trop}(\RR^2,d) \to \RR^{2n-2}
         \times \calM_{0,4,\trop}, \]
  i.e.\ $ \pi $ describes the first coordinate of the first marked point, the
  second coordinate of the second marked point, both coordinates of the other
  marked points, and the point in $ \calM_{0,4,\trop} $ defined by the first
  four marked points. Then $\deg_{\pi}(Q):= \sum_{P \in \pi^{-1}(Q)}
  \mult_{\pi}(P) $ (where $\mult_{\pi}(P)$ is defined in \cite{GM}, definition
  3.1) does not depend on $Q$ (as long as $Q$ is in general position).
\end {theorem}

\begin{proof}
  We define the map
    \[ \pi':= \ev_1^1 \times \ev_2^2 \times \ev_3 \times \cdots \times \ev_n
       \times \ft_4 : \calM_{0,n,\trop}^{\lab}(\RR^2,d) \to \RR^{2n-2} \times
       \calM_{0,4,\trop} \]
  obtained from $\pi$ by replacing $\calM_{0,n,\trop}(\RR^2,d)$ by
  $\calM_{0,n,\trop}^{\lab}(\RR^2,d)$. Then for each inverse image $P \in
  \calM_{0,n,\trop}(\RR^2,d)$ with $ \pi(P)=Q$ there exist $|\mathbb{S}_d^3|$
  different $P' \in \calM_{0,n,\trop}^{\lab}(\RR^2,d)$ with $ \pi'(P')=Q$ of
  the same multiplicity $\mult_{\pi'}(P')=\mult_{\pi}(P)$ (for the different
  labelings of the non-marked unbounded edges). Hence
    \[ |\mathbb{S}_d^3|\cdot \deg_{\pi}(Q)=\deg_{\pi'}(Q). \]
  So it is enough to show that $\deg_{\pi'}(Q)$ does not depend on $Q$. By
  proposition \ref{prop-m0nrdfan} the space $\calM_{0,n,\trop}^{\lab}(\RR^2,d)
  $ is a tropical fan. By example \ref{ex-m04} the space $ \calM_{0,4,\trop}$
  is just a tropical line in $ \RR^2 $ and thus an irreducible tropical fan by
  example \ref {ex-irred}. As $\RR^{2n-2}$ is an irreducible tropical fan
  (consisting of just one cone), too, we can conclude using proposition \ref
  {prod-irred} that $Y := \RR^{2n-2} \times \calM_{0,4,\trop}$ is an
  irreducible tropical fan. Obviously, the source and target of $ \pi' $ are
  of the same dimension. As $\pi'$ is a morphism of fans by
  propositions \ref{prop-forgetmorphism} and \ref{prop-evmorphism} we
  can apply corollary \ref{cor-image} and conclude that
    \[ \sum_{P \in |\calM_{0,n,\trop}^{\lab}(\RR^2,d)|: \pi'(P)=Q} \mult_P \pi'
       \]
  does not depend on $Q$. So it remains to show that $\mult_P \pi'=\mult_{
  \pi'}(P)$. But this follows from corollary \ref{cor-det} in the same way as
  in remark \ref {rem-indep}.
\end{proof}

\begin{remark}
  In the same way as in theorem \ref {thm-indep} the result of theorem \ref
  {thm-indep2} can of course be generalized immediately to the case of curves
  in $ \RR^r $ of arbitrary degree $ \Delta $ and with various linear
  subspaces as conditions.
\end{remark}

\end {document}